\pdfoutput=1

\documentclass[a4paper, 11pt, onecolumn, oneside]{article}
\usepackage[T1]{fontenc}
\usepackage[utf8]{inputenc}
\usepackage{makeidx}
\usepackage{graphicx}
\usepackage{float}
\usepackage{amssymb}
\usepackage{amsmath}
\usepackage{latexsym}
\usepackage{wasysym}
\usepackage{amsthm}
\usepackage{natbib}
\usepackage{}
\usepackage{amsfonts}
\usepackage{amsmath}
\usepackage{amssymb}
\usepackage{amstext}

\begin{document}

\title{ANALYSIS OF ROUND OFF ERRORS WITH REVERSIBILITY TEST AS A DYNAMICAL INDICATOR}

\author{Faranda, Davide\\
\small{\textit{Klimacampus, University of Hamburg;}}\\
\small{\textit{Grindelberg, 5, 20144 Hamburg, DE.} davide.faranda@zmaw.de}\\ \\
Mestre, Mart\'in Federico \\
\small{\textit{Grupo de Caos en Sistemas Hamiltonianos.}}\\
\small{\textit{Facultad de Ciencias Astron\'omicas y Geof\'isicas, Universidad Nacional de La Plata.}}\\
\small{\textit{Instituto de Astrof\'isica de La Plata (CCT La Plata - CONICET, UNLP).}}\\
\small{\textit{Paseo del Bosque S/N, La Plata 1900, Argentina.}}\\
\small{mmestre@fcaglp.unlp.edu.ar}\\ \\
Turchetti, Giorgio\\
\small{\textit{Department of Physics, University of Bologna.INFN-Bologna}}\\
\small{\textit{Via Irnerio 46, Bologna, 40126, Italy.} turchett@bo.infn.it}\\ \\
}

\date{}

\maketitle

\begin{abstract}
We compare the divergence of orbits and the reversibility error
for discrete time dynamical systems. These two quantities are used to explore the 
behavior of the global error induced by  round off in the computation of orbits.
The similarity of results found for any system we have analysed suggests the use of the 
reversibility error, whose computation is straightforward  since it does not require the knowledge 
of the exact orbit, as a dynamical indicator.

The statistics of fluctuations induced by round off for an ensemble of initial conditions has been 
compared with the results obtained in the case of random perturbations. 
Significant differences are observed in the case of regular orbits due to the 
correlations of round off error, whereas the results obtained for the chaotic case are nearly the same.

Both the reversibility error and the orbit divergence computed for the same number of iterations on the 
whole phase space provide an insight on the local dynamical properties with a detail comparable with other 
dynamical indicators based on variational methods such as the finite time maximum Lyapunov characteristic
 exponent, the mean exponential growth factor of nearby orbits and the smaller alignment index. 
For 2D symplectic maps the differentiation between regular and chaotic regions is well full-filled.
For 4D symplectic maps the structure of the resonance web as well as the nearby weakly chaotic regions are accurately described. 

\end{abstract}

\section{Introduction}
The discrete time dynamical systems are widely used  because 
they present a rich structure (regular, chaotic, intermittent and other
types of  orbits)  in one  or two dimensions and because the 
numerical evaluation  of the orbits is straightforward. A large number of dynamical tools known as indicators of stability allow us
to improve our understanding of dynamical systems: Lyapunov Characteristic Exponents (LCEs) are known from a long time  
\citep{wolf1985determining}, \citep{rosenstein1993practical}, \citep{skokos2010lyapunov}
as well as Return Times Statistics \citep{kac1934notion}, \citep{gao1999recurrence}, \citep{hu2004statistics} , \citep{buric2005statistics}. In the recent past many other indicators have been introduced not only to address the 
same problem of quantifying the degree of chaoticity of an orbit but also to perform the task fast. 
The Smaller Alignment Index  (SALI), widely described in \citet{skokos2002smaller} and \citet{skokos2004detecting}, 
allows to discriminate regular from chaotic orbits. Similarly, the Generalized Alignment Index (GALI), introduced
in \cite{skokos2007geometrical}, is a family of highly efficient algorithms. Besides, the Mean Exponential Growth factor of 
Nearby Orbits (MEGNO) discussed in \citet{cincotta2003phase} and \citet{gozdziewski2001global} is a quantity
that gives a fast identification of the chaoticity of the orbit while its average slope estimates 
the maximum LCE (mLCE), see \citet{2011IJNLM..46...23M} for a test of the MEGNO.
Fidelity and correlation decay are also  tools that can be  successfully used to 
characterise stability properties as explained in \citet{vaienti2007random} and \citet{turchetti2010relaxation} as well as 
Frequency Map Analysis \citep{laskar1999introduction}, \citep{robutel2001frequency}.\\ 
As it is already known, in any numerical computation of a given trajectory, there is a round off error, and it would 
be interesting to study its relation with the chaoticity of the orbit and if possible determine 
its effect with the help of some dynamical indicator.\\
The error between an exact and a numerical orbit is due to the finite
precision  used to represent real numbers  and to the arithmetics with 
round off. This is unavoidable because the length of the binary strings 
representing real numbers must not change after arithmetic operations.
The shadowing lemma is often invoked to state the existence of a true orbit
close to a numerical orbit for chaotic systems \citep{katok}, \citep{hammel1987numerical}, \citep{chow1992numerical}, \citep{chow1991numerical}. However it does not provide 
information on error growth for a given numerical orbit and the case
of regular numerical orbits is not covered by such a lemma.
The global error between the exact and the numerical orbit 
is unknown because the  first one  is not computable. Nevertheless an
estimate can be provided by  replacing the exact orbit with another one 
having  very high accuracy.  If the map is invertible the reversibility 
error can be computed  without any reference to the
exact orbit. 
 Both errors have a similar behavior, namely an average  linear growth 
for regular orbits and an average exponential growth for chaotic orbits and
consequently can be used as dynamical indicators.
Due to the correlation between the single step errors, there is a substantial 
difference with respect to the case in which the exact system is perturbed with 
random uncorrelated noise.\\
This difference has been analysed in \citet{turchetti2010relaxation}
by using the fidelity which measures the deviation of the orbits of a given map and its 
perturbation by integrating over all the initial conditions  with the appropriate measure. 
In the case of regular orbits with random perturbations the decay of fidelity is exponential 
whereas with round off errors the decay follows a power law. In the case of chaotic orbits
the asymptotic limit is approached super-exponentially for both situations (random uncorrelated
perturbations and round off).  The symplectic maps of physical interest are generally provided by the composition of a linear map with another one whose generating function is the identity plus a function of position or momentum only. In this case the
inversion is  immediately obtained in analytic form.
As most of the symplectic maps used in the literature are of this kind, 
the condition of invertibility of the map is not too restrictive. In addition, the reversibility error can also be applied to 
time-reversible systems of differential equations.\\
Here we present an  analysis of the reversibility error and its 
comparison with the divergence of orbits due to round off and other dynamical indicators of stability such as 
SALI, MEGNO and a finite time numerical estimation of the mLCE. 
Moreover we talk about the similarities and differences between an irreversibility
due to the single precision round off and one due to the application of uncorrelated random noise 
in an orbit iterated with double precision. 
For two dimensional area preserving maps the reversibility error for a 
fixed number of iterations detects the
various regions of phase space with different stability properties quite 
effectively as well as other dynamical indicators. 
Besides, the reversibility error allows to study the structure of the resonance web of a 
four dimensional symplectic map. 

%
%

\section{Round off error methods}

In a computational device a real number $x$ can be represented by a 
floating-point number $x^*$, that according to \citet{Goldberg91whatevery} and working
with base $2$, can be written as

\begin{equation}
  x^*= \pm d_0.d_1 \cdots d_{p-1} \times 2^e = \pm \sum_{k=0} ^{p-1}\; d_k 2^{e-k}  
\end{equation}
where $d_0.d_1 \cdots d_{p-1}$ is called the significand and has $p$ binary digits 
 $d_k$, whose value is $0$ or $1$, and where the exponent is an integer that
satisfies $e_{min}\le e \le e_{max}$. 
The IEEE  754 standard states that for single precision $p=24$, $e_{min}=-126$ and
 $e_{max}=127$ while for double precision $p=53$, $e_{min}=-1022$ and
 $e_{max}=1023$.

Consequently, a floating point number $x^*$ differs from the real number it represents 
and the relative error $r_p$, defined by $x^*=x(1+r_p)$, satisfies $|r_p| \leq  \epsilon \equiv 2^{-p}$, 
as analysed by \citet{Goldberg91whatevery} and \citet{knuth1973art}.
Therefore, according to IEEE 754 we have  that $\epsilon= 2^{-24}$ and $\epsilon= 2^{-53}$
for single and double precision which roughly corresponds
to 7 and 16 decimal digits, respectively. The arithmetic operations such as 
sums or multiplications 
imply a round  off, which  propagates  the error affecting each number. Round off algebraic
procedures are hardware dependent as detailed in \citet{knuth1973art}.
Unlike the case of stochastic perturbations, the error strongly depends on $x$.
Suppose we are given a map $M(x)$ then the error with respect to
the numerical map $M_*(x)$ after the first step is defined by:

\begin{equation}
  \delta_1 \equiv  \epsilon \xi_1 =x_1^*-x_1\equiv  M_*(x) -M(x) . 
\end{equation}

Analogously, we define the local error produced in the $n^{th}$ step by
\(  \delta_n \equiv \epsilon  \xi_n=M_*(x^*_{n-1}) -M(x^*_{n-1})  \)
where $\xi_n=\xi_n(x^*_{n-1})$.

The global error 

\begin{equation}
G_n= M_*^n(x)-M^n(x)
\end{equation}

accumulates all the local errors and explicit expressions can be written at
first order in $\epsilon$. In the example of a regular map we take  the translation on the torus $\mathbb{T}^1$ defined as:
\begin{equation}
  \label{torus_transl}
  M(x)= x+\omega \mod 1,
\end{equation}  
so $M_*(x)= x^*+\omega^* \mod 1$ and the global  error, which includes also the error to the mod 1 operation, becomes:

\begin{equation}
  \label{global_error}
  G_n= \epsilon\sum_{k=1}^n \;\; \xi_k =  \epsilon(n\bar \xi +w_n), 
\end{equation}

where $\epsilon\bar \xi$ is a time average defined as the limit of  $G_n/n$
for $n\to \infty$ and $w_n$ is a bounded fluctuation.

For the chaotic Bernoulli  map 
\begin{equation}
  \label{Bernoulli_map}
  M(x)=qx \mod 1,
\end{equation}
we have that $M_*(x)=qx^* \mod 1$ and that the single step and global error satisfy, 
respectively, $|\delta_1| \leq C \epsilon q$ and $|G_n| \leq C \epsilon q^n$.

In order to compute the global error we need  
the knowledge of the exact map, which is usually precluded.
A practical way to overcome this difficulty is to replace the exact map
with a map computed with an accuracy $2^{-P}$ where $P\gg p$. 
For instance, as $p=24$ corresponds to single precision  one might 
choose $P=53$ corresponding to double precision. If $p=53$ then 
one might choose $P=100$ and so forth.  There are available libraries
which allow to compute with any fixed number of bits or significant
decimal digits. The computation of the reference orbit is expensive
if high precision is used, but there is no other way to evaluate the global error. As a consequence the ``exact" orbit is achievable 
for a definite number of iterations which depends on $P$ and the nature of 
the map. Taking into account what we have just mentioned, in the forthcoming numerical 
experiments, we will use the divergence of orbits, defined by:
\begin{equation}
  \Delta_n= M_S^n(x)-M_D^n(x),
   \label{eq:def_diver}
\end{equation}
where $M_S$ and $M_S$ stand for single and double precision iterations respectively.

If the map is invertible there is another option to overcome the difficulty of 
not possessing the true map.  We define the reversibility error as

\begin{equation}
  R_n= M_*^{-n}\circ M_*^n (x)-x
  \label{eq:def_rever}
\end{equation}

which is non zero since the numerical inverse $M_*^{-1}$  of the map 
is not exactly the inverse of $M_*$ namely $M_*^{-1}\circ M_*(x)\not = x$.
Obviously the reversibility error is much easier to compute than the divergence 
of orbits (if we know explicitly the inverse map) and the information it provides 
is basically the same as the latter. 
Both quantities give an average linear growth for a regular map together with an 
exponential growth for a chaotic map having positive Lyapounov exponents and strong mixing properties. When computing $R_n$  we will set $M_* = M_S$
in order to compare with $\Delta_n$.\\

\section{Variational methods}

In the forthcoming sections we will compare the performance of the indicators presented above with
three well known and widely accepted dynamical indicators that are based on the behavior
of the solution of the variational equations of the system. These are the finite time mLCE,
the cumulative moving time average of MEGNO, and SALI.

Let us briefly state them for discrete time dynamical systems of the form:
\begin{equation}
\mathbf{x}_{n+1}= \mathbf{f}(\mathbf{x}_n),
\end{equation}
where $\mathbf{x}_n$ is the state vector at time $n$ and $\mathbf{f}$ is a vector valued function.

The concomitant discrete time variational equations, also called tangent map dynamics, associated
to a given orbit $\{\mathbf{x}_n\}_{n\in \mathbb{N}}$ are the following:
\begin{equation}
\mathbf{v}_{n+1}= \mathbf{Df}(\mathbf{x}_n)\cdot \mathbf{v}_n,
\end{equation}
where $\mathbf{Df}(\mathbf{x})$ is the Jacobian matrix of the function $\mathbf{f}$ and 
$\mathbf{v}_n$ is a deviation vector at time $n$.

\citet{skokos2010lyapunov} makes a historical review of the definition of the LCEs and its
connection with the divergence of nearby orbits. He also states the theorems that guarantee
the existence of the spectrum of LCEs and, in particular, the current definition of the mLCE
in  terms of the solution of the variational equations: 

\begin{equation}
\label{mLCE}
\mbox{mLCE}\equiv \lim_{n \to \infty}  \frac{1}{n}\ln \frac{||\mathbf{v}_n||} {||\mathbf{v}_0||},
\end{equation}

with $||.||$ some norm. For a chaotic orbit the mLCE is positive and this implies 
an exponential divergence of nearby orbits. On the other hand, for regular orbits mLCE is 
zero. 
 
In order to have a numerically computable quantity we define the finite time mLCE 
at time $n$ as 
\begin{equation}
   \label{eq:def_mLCE}
   \mbox{mLCE}(n)\equiv \frac{1}{n} \ln \frac{||\mathbf{v}_n||} {||\mathbf{v}_0||} = 
   \frac{1}{n} \sum_{k=1}^n \ln \frac{||\mathbf{v}_k||} {||\mathbf{v}_{k-1}||} ,
\end{equation}
 so that equation  Eq. ~(\ref{mLCE}) can be reformulated as
\begin{equation}
\mbox{mLCE} = \lim_{n \to \infty} \mbox{mLCE}(n).
\end{equation}
 
\citet{cincotta2003phase} defines a biparametric family of MEGNO indicators:

\begin{equation}
\label{MEGNO}
Y_{m,j}(n) = (m+1)n^j \sum_{k=1}^n k^m \ln \frac{||\mathbf{v}_k||} {||\mathbf{v}_{k-1}||} ,
\end{equation}
where $m$ and $j$ are integer numbers. They made experiments with $Y_{2,0}$, $Y_{3,1}$ and $Y_{1,-1}$
and concluded that the last one allows both a fast classification between chaotic and regular orbits, 
and a clear identification of stable and unstable periodic orbits. Due to this fact, we will use 
$Y \equiv Y_{1,-1}$ throughout the rest of this article. In order to reduce the fast oscillations
that the time evolution of the MEGNO presents, in \cite{cincotta2003phase} they use a time average
of this quantity, namely:

\begin{equation}
  \label{eq:def_av_MEGNO}
  \bar{Y}(n) = \frac{1}{n} \sum_{k=1}^n Y(k).
\end{equation}

Theoretically, the asymptotic evolution of $\bar{Y}(n)$ for any dynamical regime can be put into a single expression:
\begin{equation}
  \label{theory_MEGNO}
  \bar{Y}(n) \approx \frac{\mbox{mLCE}}{2} n + c,
\end{equation}
where $c\approx 0$ and $2$ for chaotic and regular motion, respectively.

We also compare our results with the SALI, introduced by \citet{skokos2001alignment},
that measures the degree in which a pair of initially linearly independent deviation vectors tend to 
become aligned. The underlying principle is that, for a chaotic orbit, a deviation vector under the
tangent map dynamics changes in order to become aligned  with the instantaneous most unstable direction. 
In other words, for almost every pair of initial deviation vectors ($\mathbf{v}_0$, $\mathbf{u}_0$), 
the more chaotic the orbit the faster that the angle between them will reduce to zero. 
In the case of regular orbits, the behavior depends strongly on
the dimensionality of the map: in maps with dimension $\geq4$ the deviation vectors  generally remain unaligned 
so the SALI tends to a positive non-zero value, while in 2D maps the two deviation vectors tend to align 
with a time rate that follows a power law. See \citep{skokos2004detecting} and  
\citep{skokos2007geometrical} for numerical tests of the SALI and its generalisation, the GALI family.
Denoting  the Euclidean norm with $||.||$ the SALI is defined as:

\begin{equation}
   \label{eq:def_sali}
  \mbox{SALI}(n)= \min \left\{  \left\| \frac{\mathbf{v_n}}{\| \mathbf{v_n} \|} + \frac{\mathbf{u_n}}{\| \mathbf{u_n} \|}   \right\|, \left\| \frac{\mathbf{v_n}}{\| \mathbf{v_n} \|} -\frac{\mathbf{u_n}}{\| \mathbf{u_n} \|}  \right\| \right\}
\end{equation}

To control the exponential increase of the norm of the vectors and avoid overflow problem, 
\citet{skokos2004detecting} have normalised them, at every time step, keeping their norm equal to 1. 
Defined like this SALI$(n) \in [0,\sqrt{2}]$ and SALI = 0 if and only if the two normalised vectors have 
the same direction, being equal or opposite.

%
%

\section{Error behavior on simple  maps}

From the linear and exponential bounds on the errors 
it follows that the numerical orbit remains close to the exact one 
for a number of iterations proportional to $1/\epsilon$ in the regular 
case and to $\ln(1/\epsilon)$ in the chaotic one.  

We consider two types of models where the error grows linearly and exponentially
respectively. The first one is the translation on the torus $\mathbb{T}^1$
defined by Eq. ~(\ref{torus_transl}).

This is equivalent to the rotations on the unit circle defined by 
the map

\begin{equation}
  \label{torus_rotat}
  R \mathbf{x} =
  \begin{pmatrix}
    \cos(2\pi\omega)  &   -\sin(2\pi\omega) \\
    \sin(2\pi\omega)  &   \cos(2\pi\omega) 
  \end{pmatrix} 
  \begin{pmatrix}
    \cos(2\pi x) \\
    \sin(2\pi x) 
  \end{pmatrix}.
\end{equation}

The correspondence between the sequences $x_n=M(x_{n-1})$ and
$\mathbf{x} _n= R \mathbf{x}_{n-1}$ is evident after writing
 $\mathbf{x}_n=(\cos(2\pi x_n), \sin(2\pi x_n))$.
In spite of the rigorous mathematical equivalence between the translations
on the torus and the rotations on the unit circle, there is not 
such similitude in the numerical evaluation of these maps, as will be explained in 
the following paragraph.
 
\begin{figure}[ht!]
  \centering  
  \includegraphics[width=0.45\textwidth]{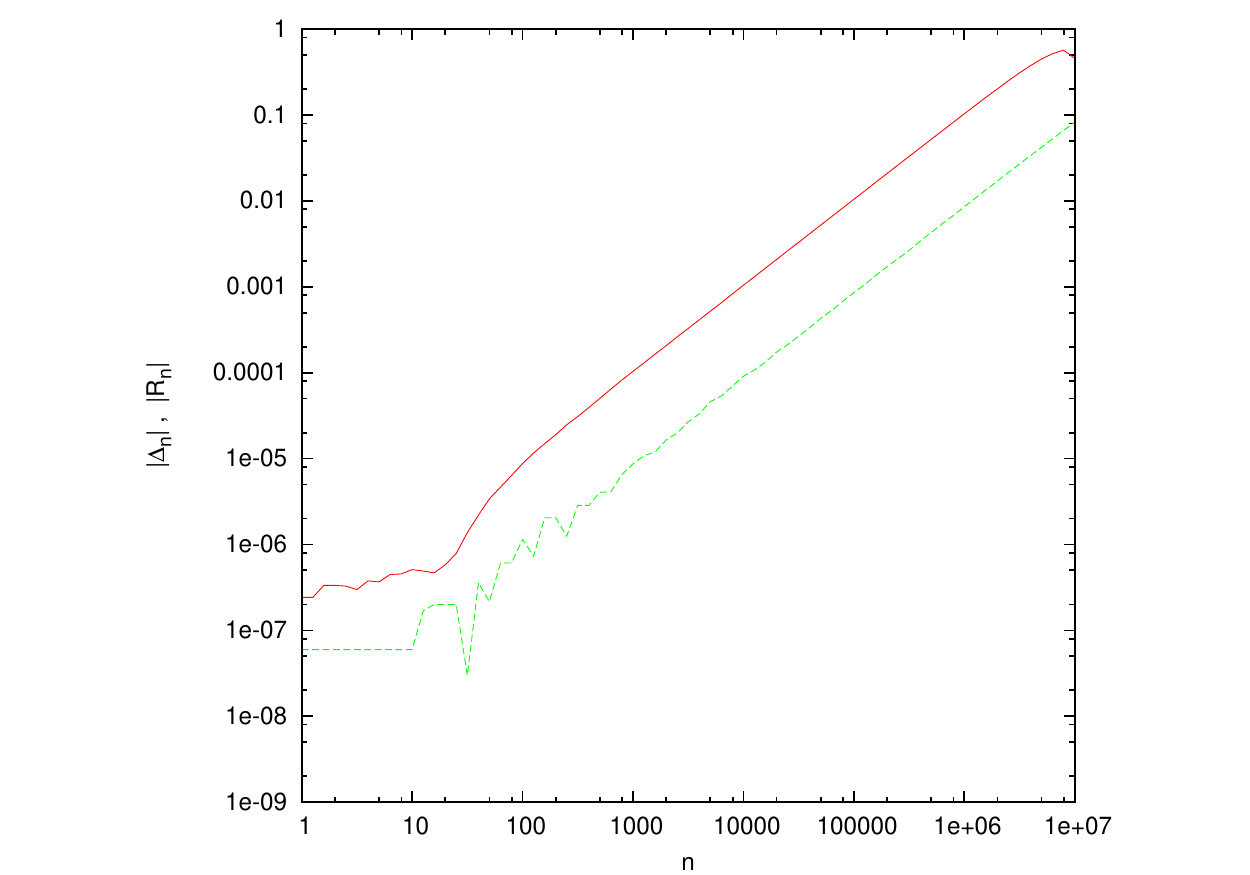}
  \hspace{1cm}
  \includegraphics[width=0.45\textwidth]{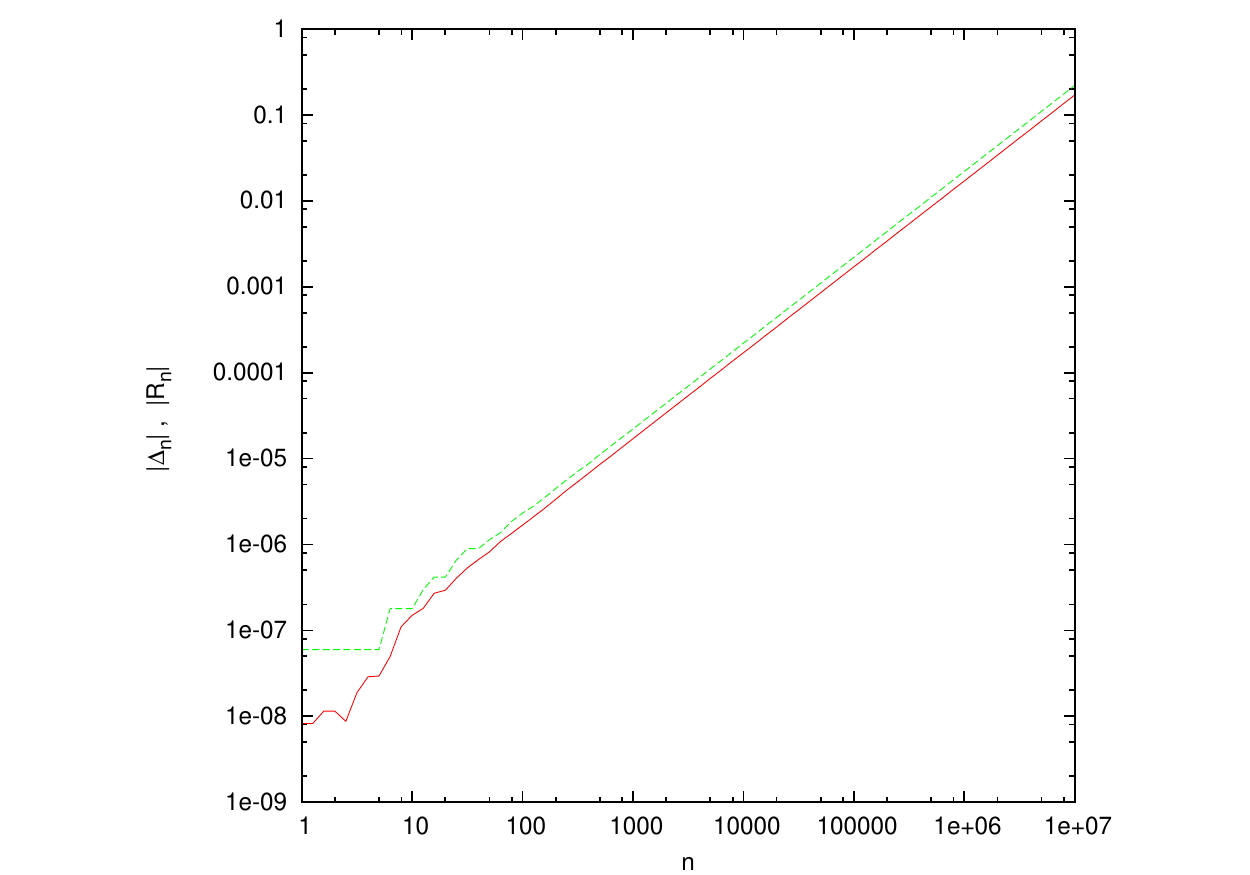}
  \caption{Comparison between divergence (red) and reversibility (green) for both torus 
    translations (left) and rotations (right), using a frequency $\omega=\sqrt{2}-1$ and 
    initial condition $x_0=0.7$. In all four cases the average slope in log-log scale is nearly one.}
  \label{fig: comparison_SP}
\end{figure}

For map (\ref{torus_transl})  the divergence of orbits and the reversibility 
error have a linear growth, basically due to the the fact that $\omega^* \neq \omega$, 
as shown theoretically for the global error in Eq.~(\ref{global_error}). However there 
are architectures and/or compilers in which the reversibility error may be zero for this simple map
that only requires the computation of sums and modulus operations. 
Without knowing the precise implementation of these peculiar operations  it is impossible to 
establish a priori with which compiler and in which architecture we may obtain this result. 
We believe this anomaly is due to the peculiar arithmetic operations involved
and  that  it does not occur for a generic map.  To support this claim we
have checked that the reversibility error never vanishes for the map (\ref{torus_rotat}) 
which involves multiplications and the evaluation  of trigonometric
functions.\\
In Figure~\ref{fig: comparison_SP}-left we compare the divergence and reversibility error for 
the torus translation with $\omega=\sqrt{2}-1$ and $x_0=0.7$. 
We can see that both errors approximately satisfy the same expression: in log-log scale it is close to a straight 
line with unitary slope. This implies that both quantities are linear functions of time.
We have checked that the behavior is similar for almost every initial condition and frequency but in some cases the fluctuations around the average linear growth have a larger amplitude. 
Figure~\ref{fig: comparison_SP}-right shows that for the map (\ref{torus_rotat}) there is an analogous
behavior. However, in this case the straight line that corresponds to orbit divergence is a bit shifted 
upwards with respect to the one of reversibility error. This is telling us that the orbit divergence's linear
growth has a bigger slope.

   The fluctuations $w_n$ have variance $<w_n^2>$ computed 
respect to an ensemble of initial conditions $x_0$  which first  grows linearly, 
but rapidly  saturates to a constant value with small oscillations
( see Figure~7 in \citep{turchetti2010relaxation} ).
In this respect it is quite different with respect to a random 
perturbation $\epsilon \xi_n$ where $\xi_n$ are independent variables,
since in this case we  have $w_n=\xi_1+\ldots+\xi_n$ and the variance is 
given by $<w_n^2> = n <\xi^2>$  for any value of $n$.  This different behavior
is reflected in the decay of fidelity \citep{vaienti2007random} which follows a power law for
round off errors and an exponential law for random perturbations \citep{turchetti2010asymptotic} , \citep{turchetti2010relaxation}.
This shows that the round off errors decorrelate very slowly unlike the random 
errors which are uncorrelated.

For chaotic maps, characterised by an  exponential increase in the distance 
between two nearby orbits,  the reversibility error due to round off has the same growth.
Bounded two dimensional chaotic maps, like the cat map, exhibit this behaviour during a short 
time before reaching saturation.

\section{The case of the standard map}
As an example of generic 2D map we have chosen the standard map in a torus:

\begin{equation} 
  y_{n+1}=y_n+ \lambda \sin(x_n) \mod 2\pi; 
  \qquad x_{n+1}=x_n+y_{n+1}\mod 2\pi. 
  \label{standard_map}
\end{equation}

For very low values of $\lambda$ the divergence of orbits
has an average growth linear with $n$ as  for the translations on
the torus and they depend weakly on the initial condition. 
As $\lambda$ is increased a power law $n^\alpha$ with $\alpha>1$ 
is observed and the dependence on the initial conditions becomes appreciable.
Reversibility error also shows this behavior for $\lambda<<1$.

\begin{figure}[H]  
  \centering  
  \includegraphics[width=0.45\textwidth]{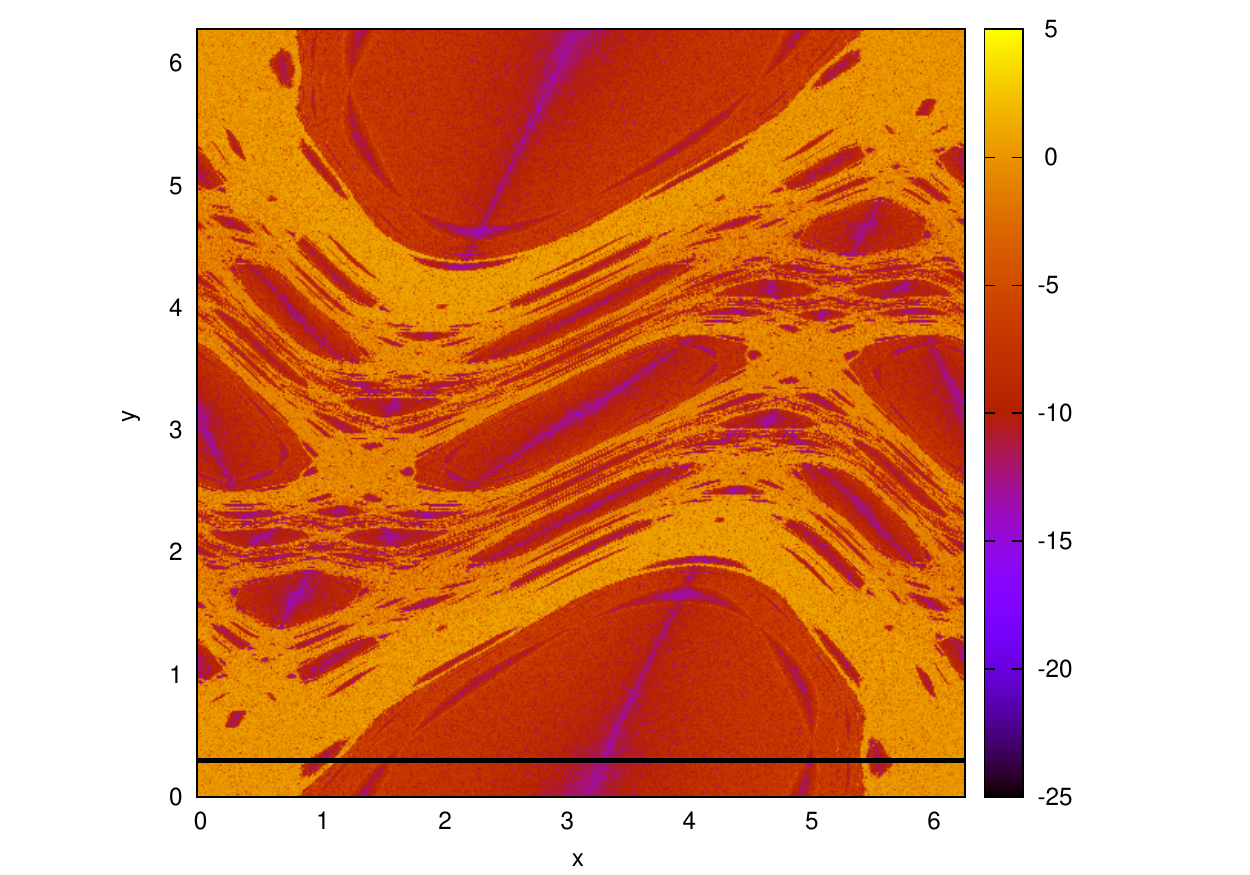} 
  \hspace{1cm}
  \includegraphics[width=0.45\textwidth]{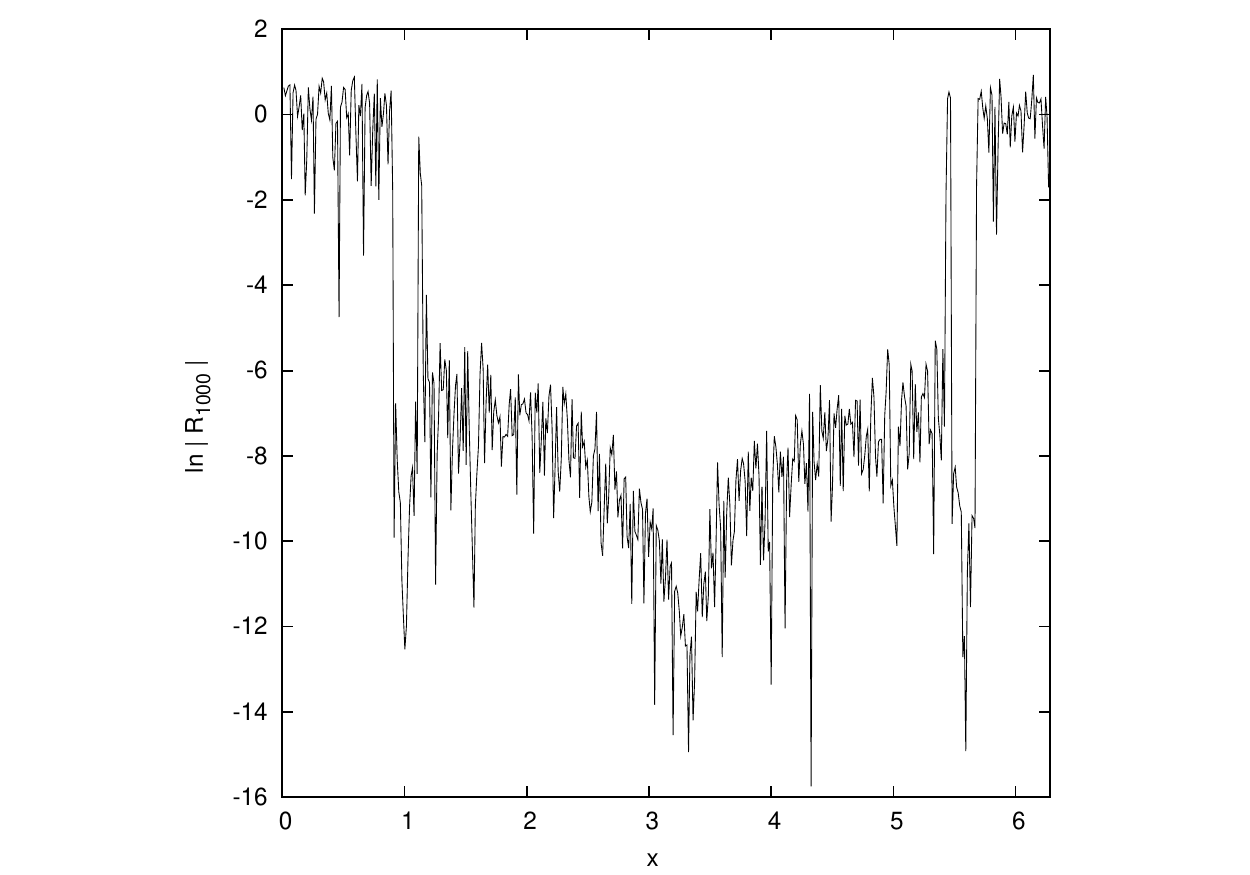}
  \caption{ $\ln(R_{10^3})$ displaying the dynamical structure of the phase space of the Standard Map (left) and for torus section  $y=0.3$ (right).}
  \label{fig: rever_SM}
\end{figure}

\begin{figure}[H]  
  \centering  
  \includegraphics[width=0.45\textwidth]{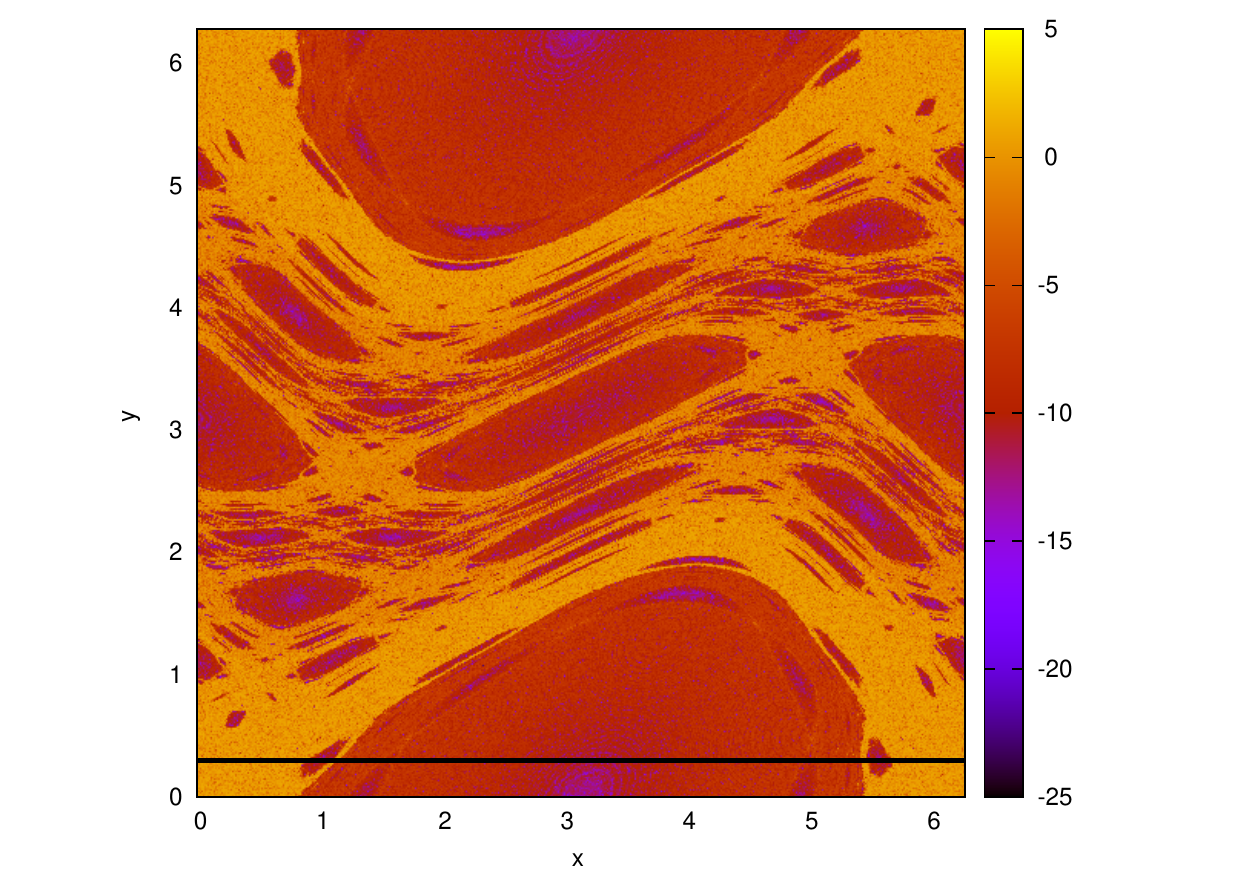} 
  \hspace{1cm}
  \includegraphics[width=0.45\textwidth]{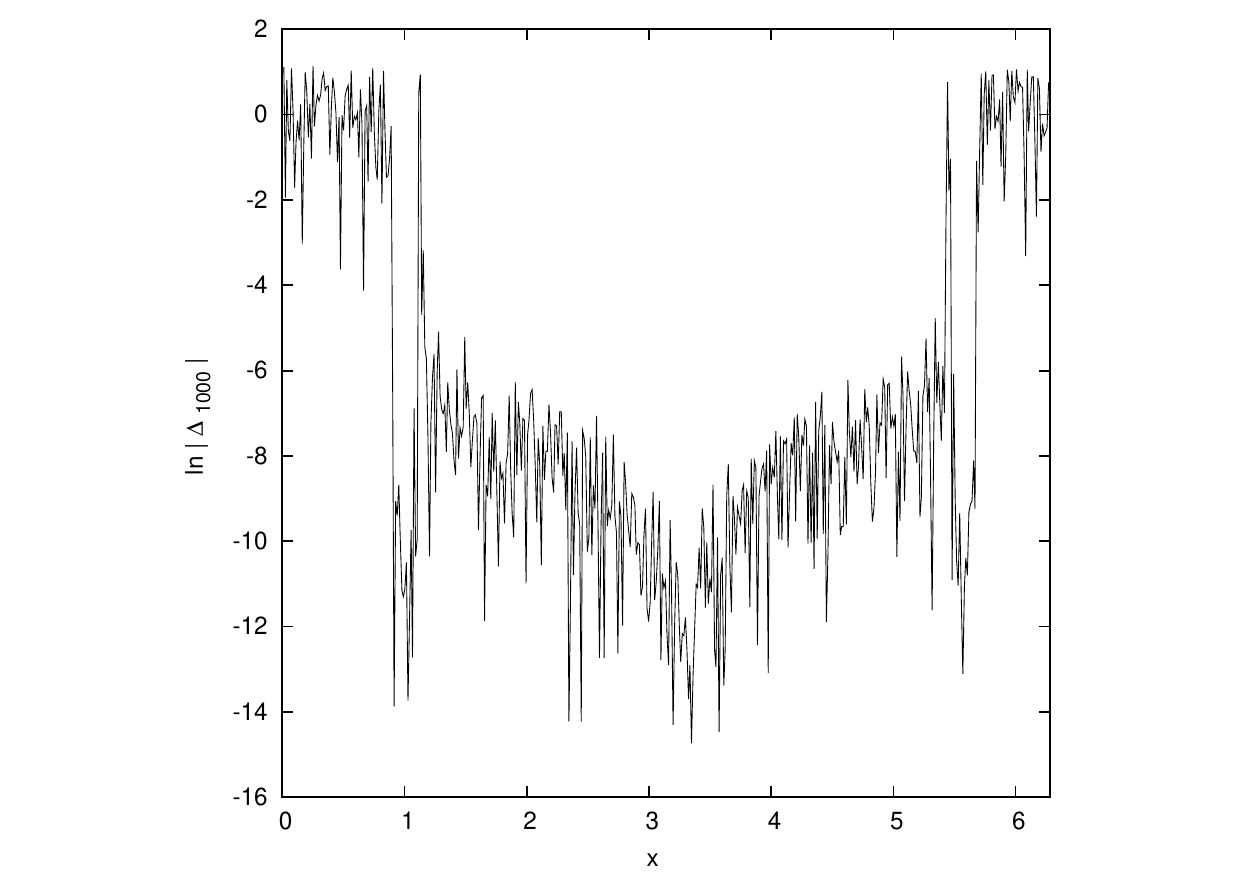}
  \caption{ $\ln(\Delta_{10^3})$ displaying the dynamical structure of the phase space of the Standard Map (left) and for torus section  $y=0.3$ (right).}
  \label{fig: diver_SM}
\end{figure}

\begin{figure}[H]  
  \centering  
  \includegraphics[width=0.45\textwidth]{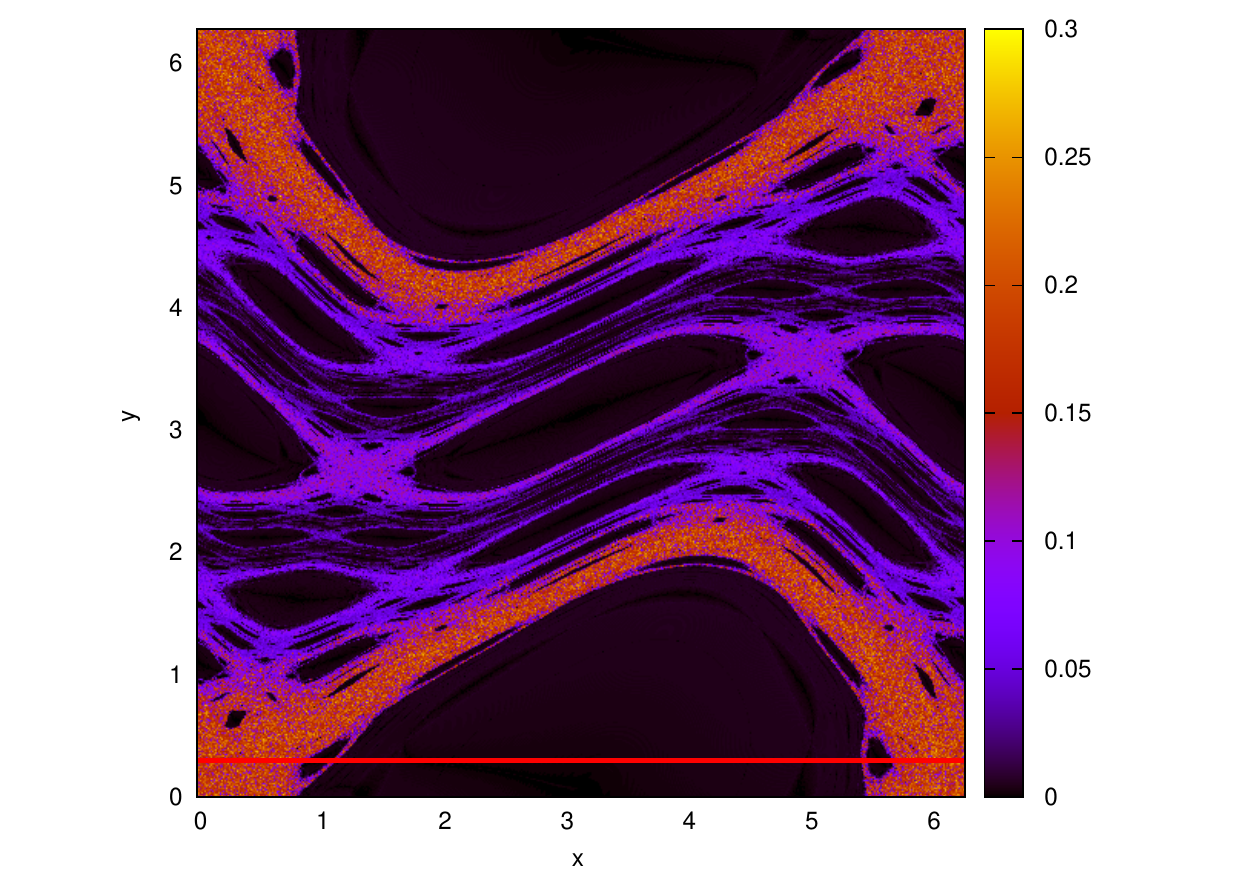} 
  \hspace{1cm}
  \includegraphics[width=0.45\textwidth]{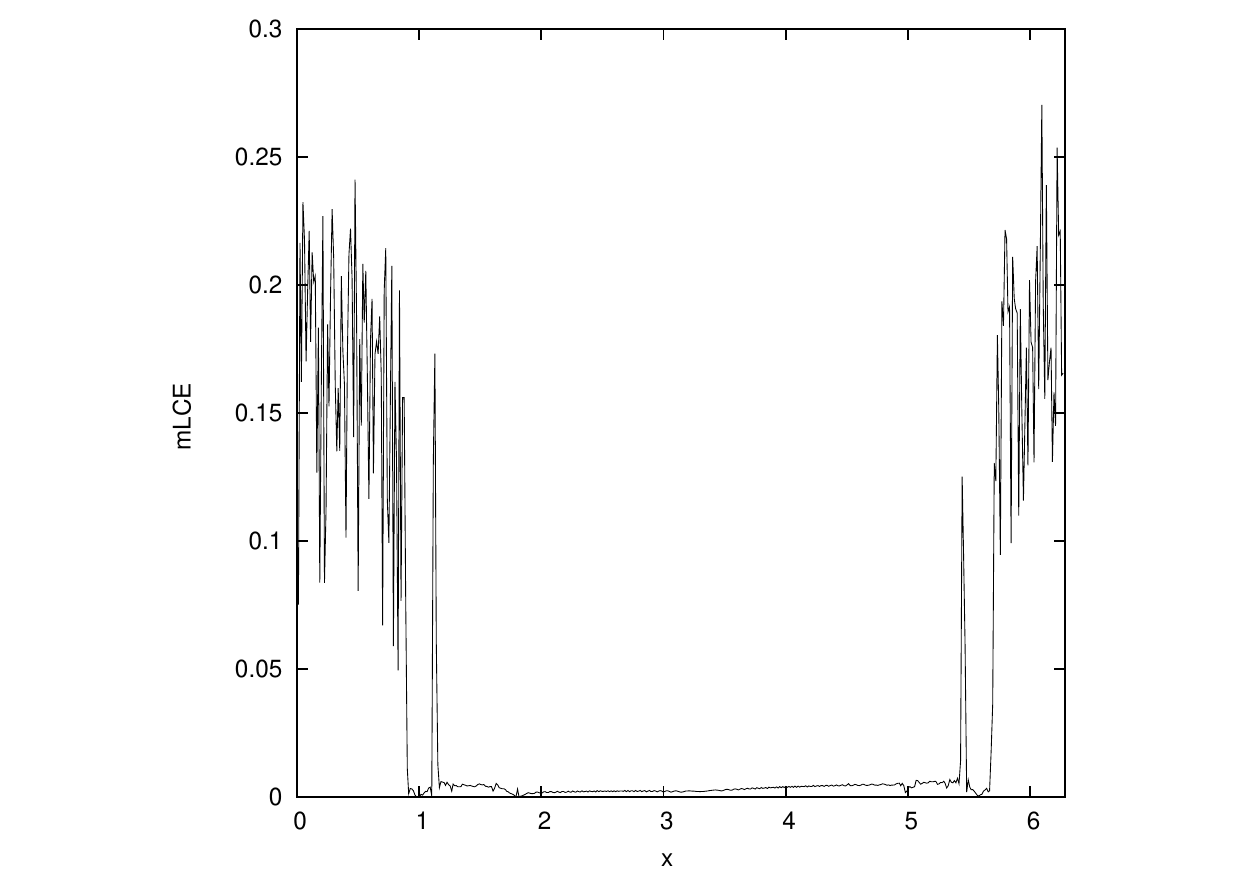}
  \caption{ $\mbox{mLCE}(10^3)$, displaying the dynamical structure of the phase space of the  Standard Map (left) and for torus section  $y=0.3$ (right).}
  \label{fig: mLCE_SM}
\end{figure}

\begin{figure}[H]  
  \centering  
  \includegraphics[width=0.45\textwidth]{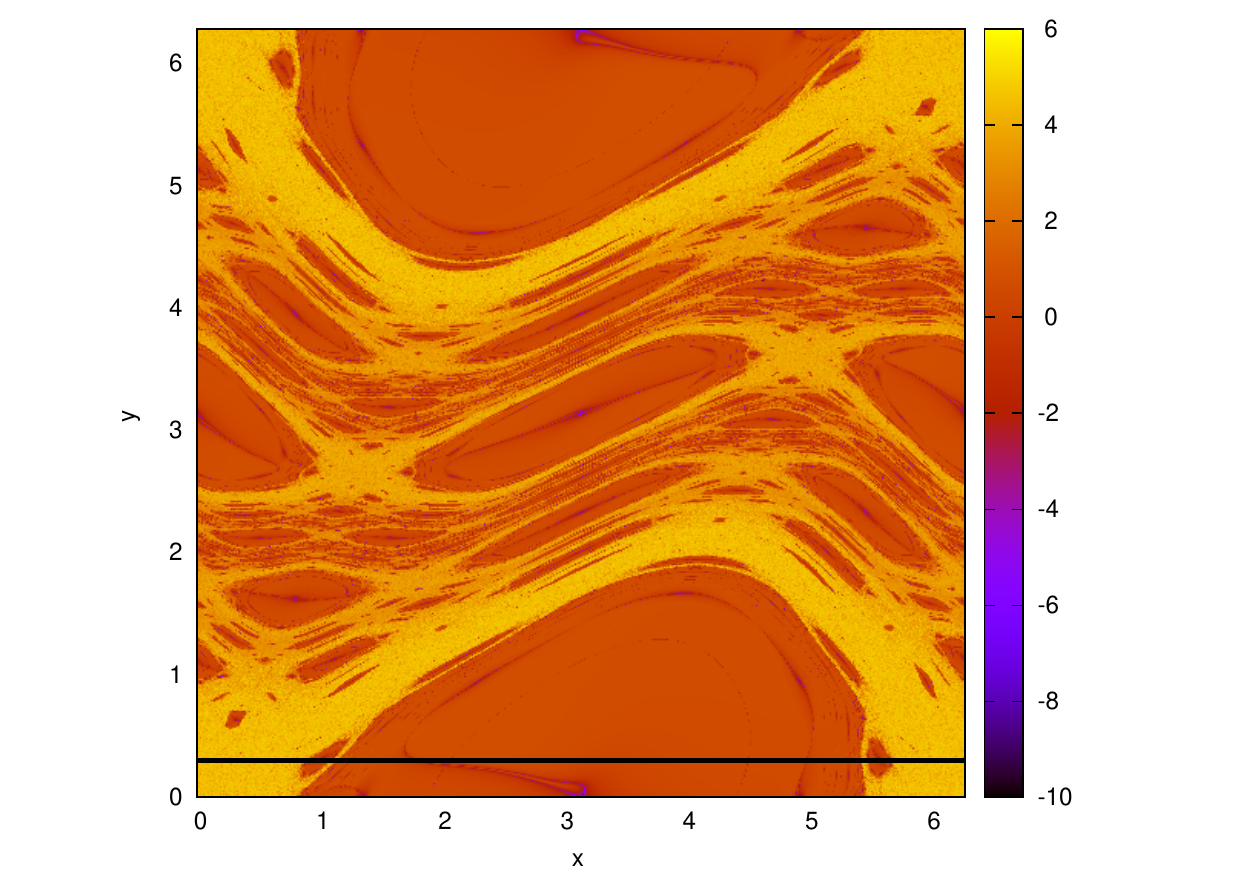} 
  \hspace{1cm}
  \includegraphics[width=0.45\textwidth]{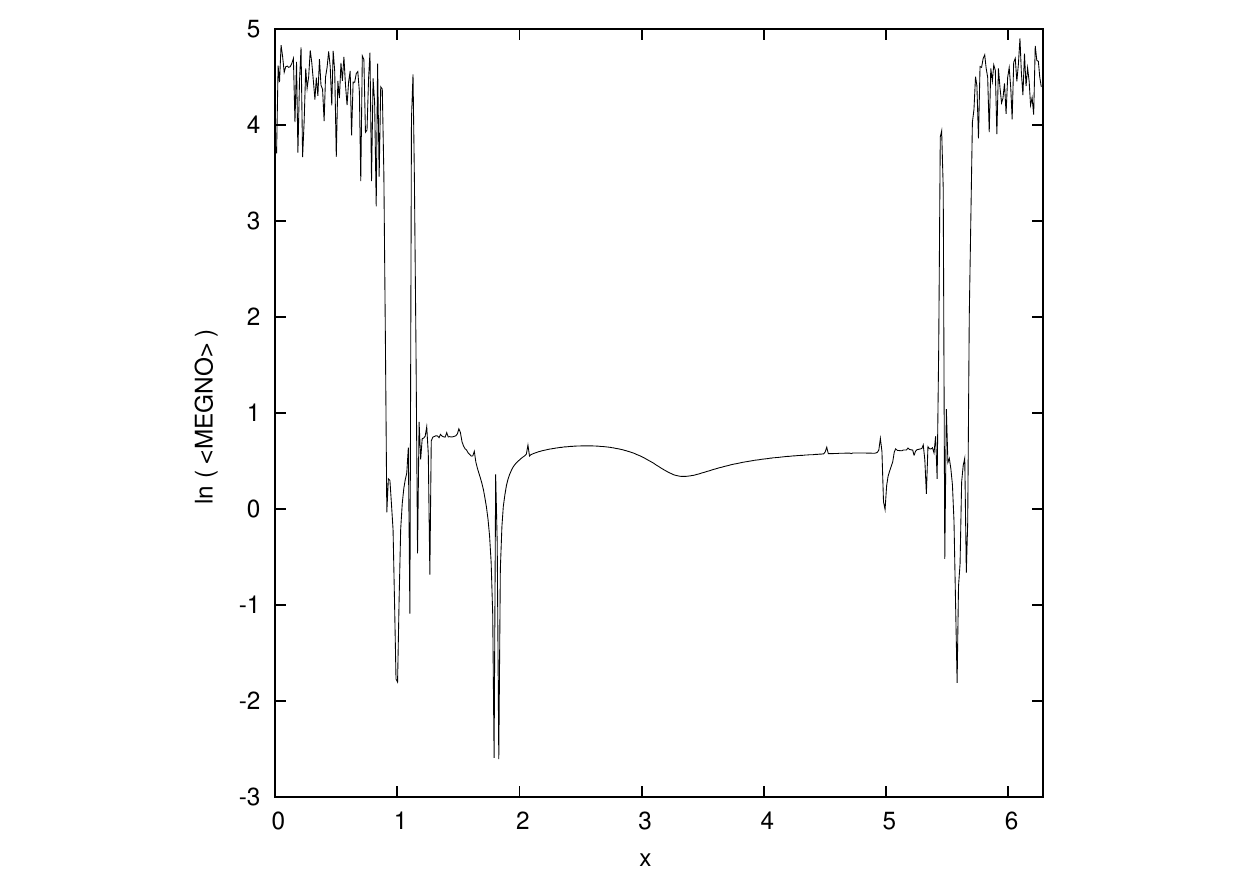}
  \caption{ $\ln(\bar{Y}(10^3))$, displaying the dynamical structure of the phase space of the 
 Standard Map (left) and for torus section  $y=0.3$ (right).}
  \label{fig: av_MEGNO_SM}
\end{figure}

\begin{figure}[H]  
  \centering  
  \includegraphics[width=0.45\textwidth]{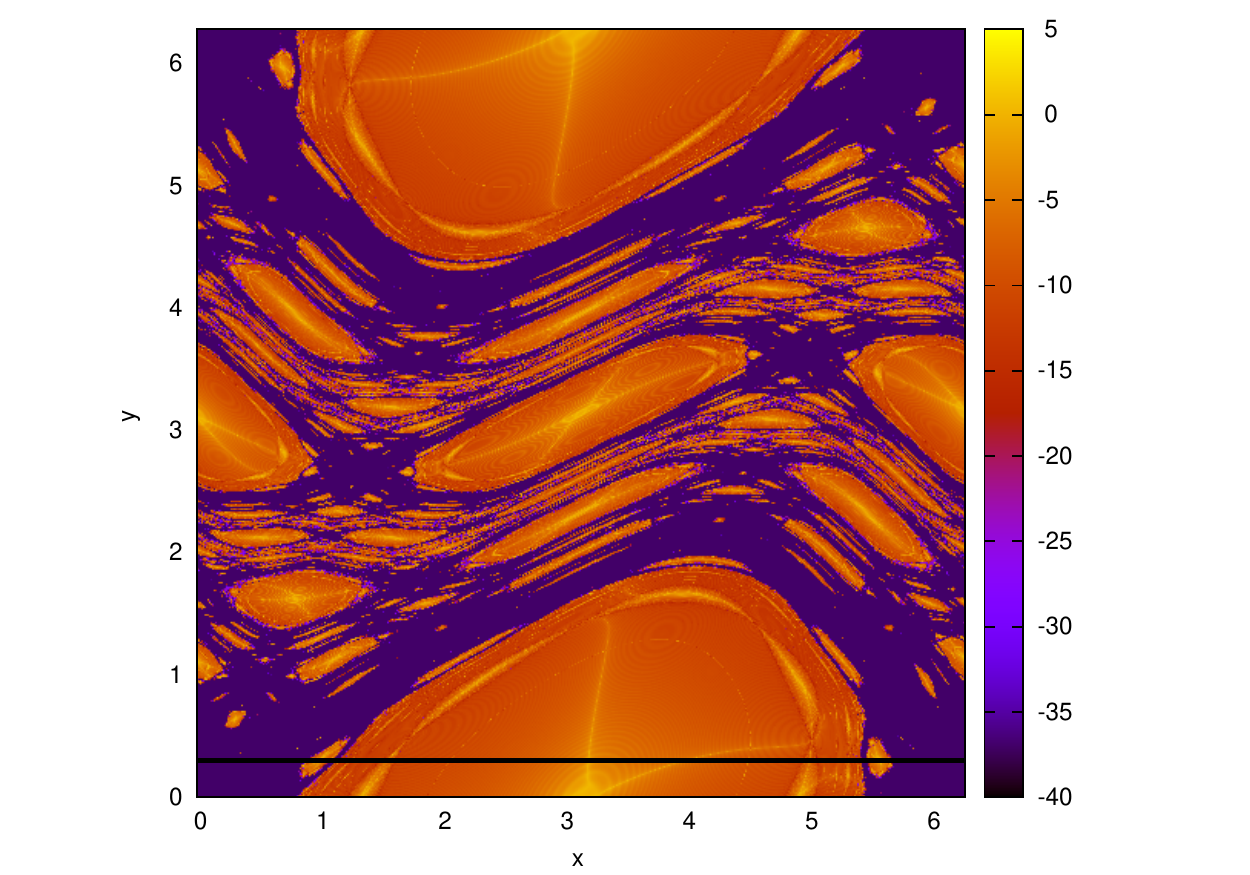} 
  \hspace{1cm}
  \includegraphics[width=0.45\textwidth]{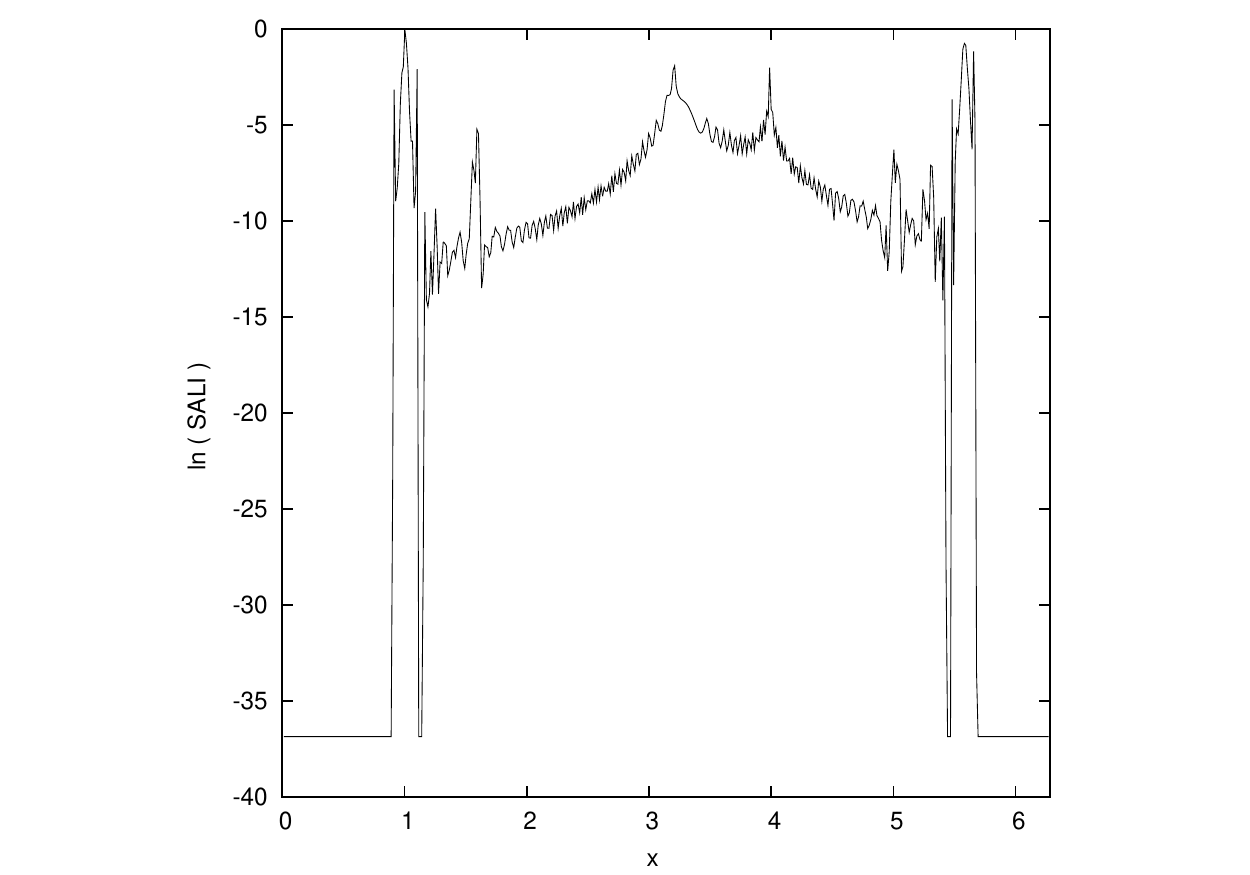}
  \caption{ $\ln($SALI($10^3$)) displaying the dynamical structure of the phase space of the Standard Map (left) and for torus section  $y=0.3$ (right).}
  \label{fig: sali_SM}
\end{figure}

For $\lambda $ approaching one we have coexistence of regular and chaotic regions so that
the domain is splitted into several islands of stable orbits and a chaotic sea.  
In figures ~\ref{fig: rever_SM}, ~\ref{fig: diver_SM}, 
~\ref{fig: mLCE_SM}, ~\ref{fig: av_MEGNO_SM}
and ~\ref{fig: sali_SM} we show, respectively, the value of the 
reversibility error~(\ref{eq:def_rever}), the orbit divergence~(\ref{eq:def_diver}), 
the finite time mLCE~(\ref{eq:def_mLCE}), the time average of the MEGNO~(\ref{eq:def_av_MEGNO})  
and the SALI~(\ref{eq:def_sali}) obtained iterating eq.~(\ref{standard_map}).

The standard map with $\lambda=0.971635$ has been iterated  $n=10^3$ times to compute all the dynamical indicators.
Even if $n<10^3$ might be used to describe the chaotic region, the chosen value $n=10^3$ allows to highlight differences within regular
regions where the growth is very slow.

The pictures in the left side show the value of  the dynamical indicators in a chromatic scale  for a grid of $500 \times 500$ 
initial conditions in the two dimensional torus whereas the right pictures show the corresponding graph for 
a fixed value of the action
variable ($y=0.3$) corresponding to the horizontal line in the left figures. Except for mLCE  we have used the natural  logarithm 
of the absolute value for all the dynamical indicators.

Figure~\ref{fig: rever_SM} and Figure~\ref{fig: diver_SM} were obtained by taking into account the error only in 
the action variable. It is evident that both the reversibility error and the orbits divergence (with respect to the exact one) 
discriminate regular from chaotic orbits. Minor differences exist within the islands of stability: 
the orbit divergence approaches the minimum value close to the centre, whereas the minimum of the reversibility error 
occurs on a line crossing the islands possibly due to a mechanism similar to the one reported by \citet{barrio2009spurious} which analysed spurious errors for variational methods.

The plots obtained for mLCE and the MEGNO (figures~\ref{fig: mLCE_SM} and~\ref{fig: av_MEGNO_SM} respectively) 
show no structure within the resonance islands whereas variations and fluctuations appear in the chaotic
regions as for the previous indicators.

Figure~\ref{fig: sali_SM} shows the logarithm of the SALI value. 
Due to the fact that this indicator converges to zero extremely fast for chaotic orbits,
we have used a cut off value at SALI $= 10^{-16}$.  The presence of this cut-off erases any structure within 
the strongly chaotic region.

To summarise, all the indicators discriminate regular and chaotic regions but their sensitivity in these regions is different .
 
Another aspect to consider when comparing the efficiency of chaos indicators is the computational cost.
Each variational method needs to iterate both the map and the tangent map forward for $n$ steps. 
The tangent map is the computationally most expensive  since it needs 
the evaluation of the Jacobian matrix at every step. When computing SALI, two deviations vectors must be simultaneously computed.
We also use two deviation vectors in evaluating MEGNO, selecting at every step the one that stretches more, in order to reduce 
the probability of having a vector almost orthogonal to
the most unstable direction. In the case of the computation
of the mLCE we have used the orthogonalization algorithm developed by~\citet{benettin1980lyapunov}.
On the other hand, we notice that the reversibility error method requires only the iteration of the map whereas the orbit divergence method
requires  the iteration  of the single and double precision (or double and higher precision) map. As a consequence, the computationally most economic method is the one based on the reversibility error which does not require any algorithm except the evaluation of the map itself.
Typically, the relevant information can be extracted from a computation in single precision.

\begin{figure}[H]  
  \centering  
  \includegraphics[width=0.45\textwidth]{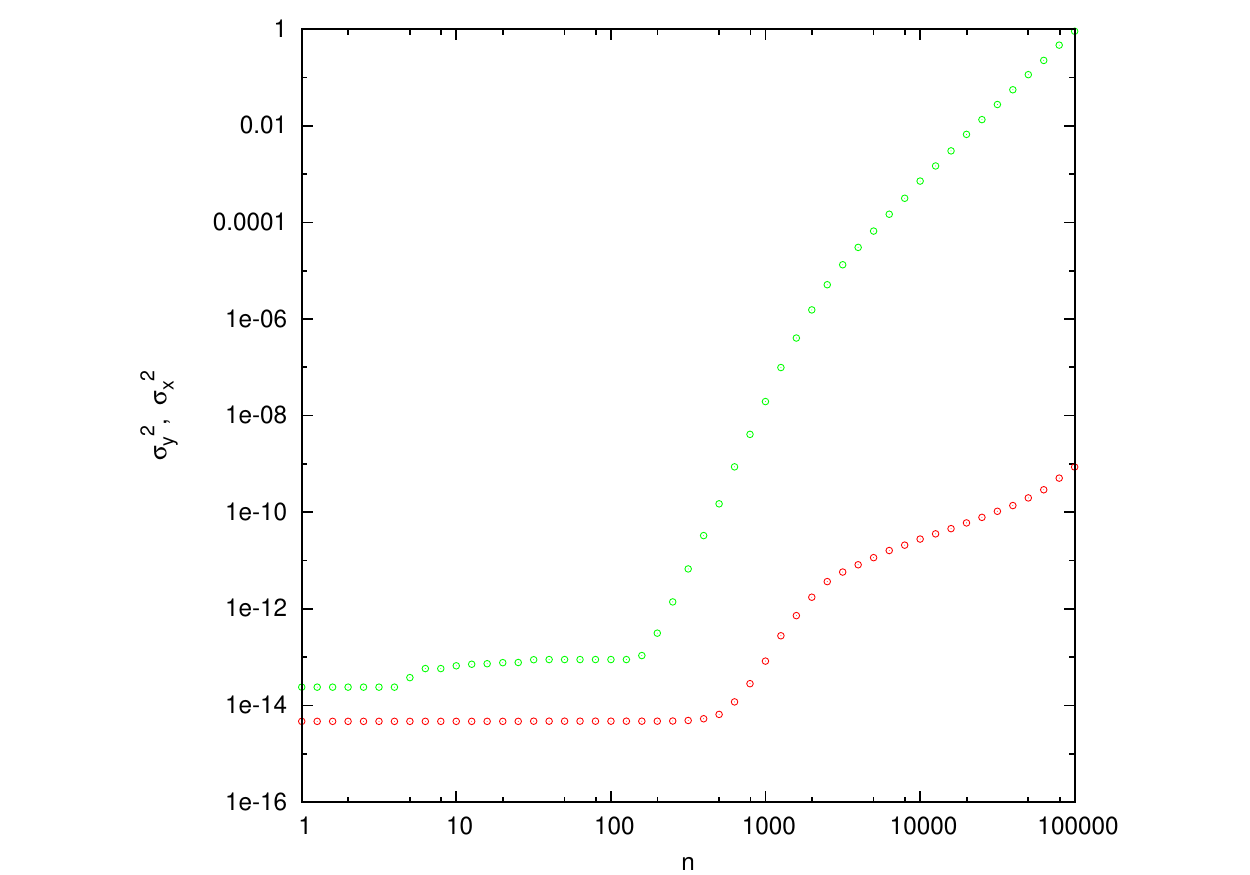} 
  \hspace{1cm}
  \includegraphics[width=0.45\textwidth]{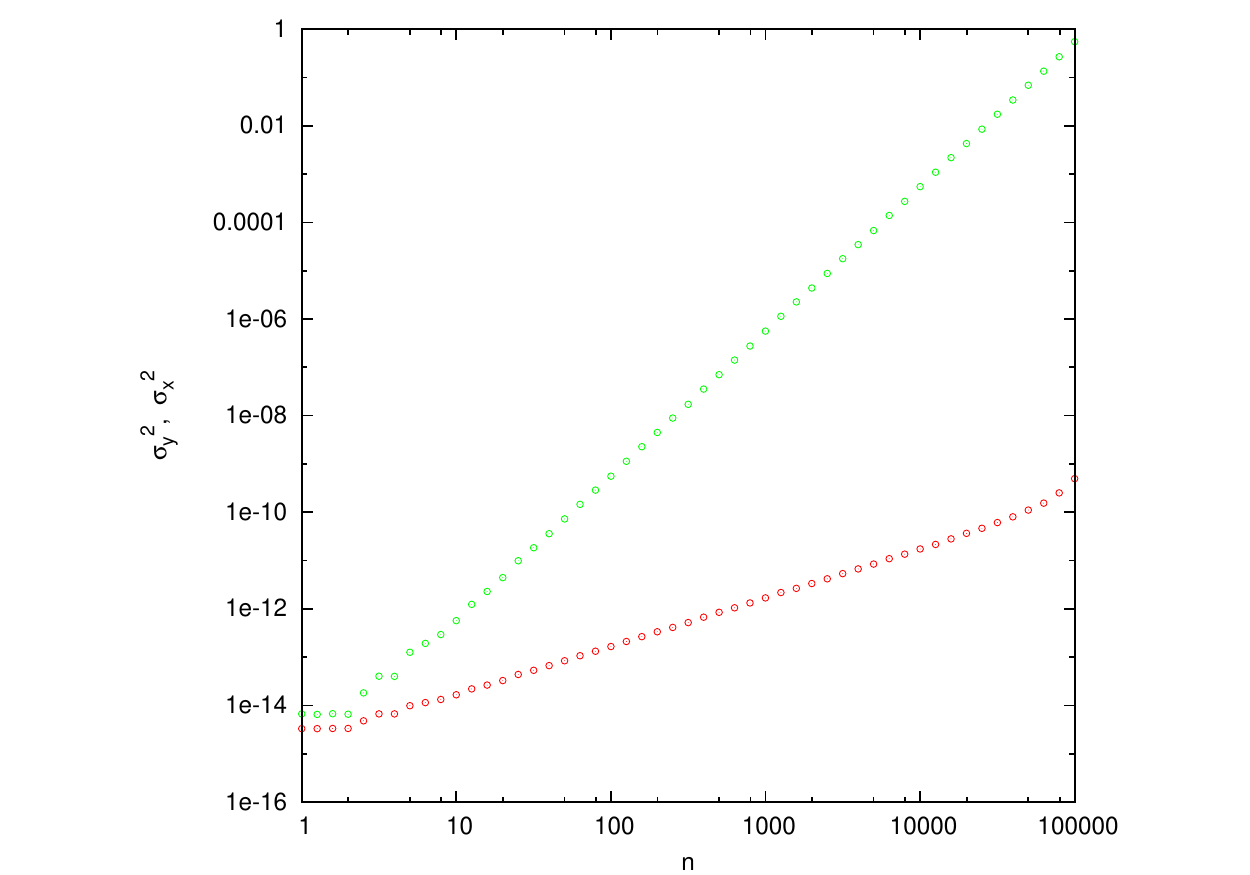}
  \caption{Evolution of the variance of a distribution of reversibility errors in action (red) and angle (green) 
    variables for the standard map with $\lambda=10^{-4}$. Left: round off error. Right: stochastic uncorrelated perturbations of 
    amplitude $10^{-7}$.}
  \label{fig: rever_variances}
\end{figure}

From now on, within this section, we will focus on some ensemble quantities in order to compare the effect
of round off errors with the effect of random perturbations in the standard map. 
For the standard map with $\lambda=10^{-4}$ all the orbits are regular and we follow
the evolution of an ensemble of 10001 initial conditions 
randomly chosen in $(x,y)\in[1.5,~1.5 + 10^{-3}]\times[\pi,~\pi+10^{-3}]$. For each iteration we compute
the variance of $R_{n}$ in action (${\sigma_y}^2$) 
and angle (${\sigma_x}^2$) variables (see in Figure~\ref{fig: rever_variances}-left) and compare it against
the same quantities 
of a double precision orbit stochastically perturbed with
a uniform uncorrelated noise of amplitude $10^{-7}$ (shown in
Figure~\ref{fig: rever_variances}-right). In the latter case we have found that ${\sigma_y}^2$ and 
${\sigma_x}^2$ grow according a power law  with exponents one and three respectively up to numerical uncertainties.

It is interesting to compare the two pictures in Figure~\ref{fig: rever_variances} with each other. 
For $n$ large we see that  the behaviors of variances due to round off and random perturbations are similar. 
The presence of a transient for the reversibility case is very likely due to the 
initial presence of correlation.

\begin{figure}[H]  
  \centering
  \includegraphics[width=0.45\textwidth]{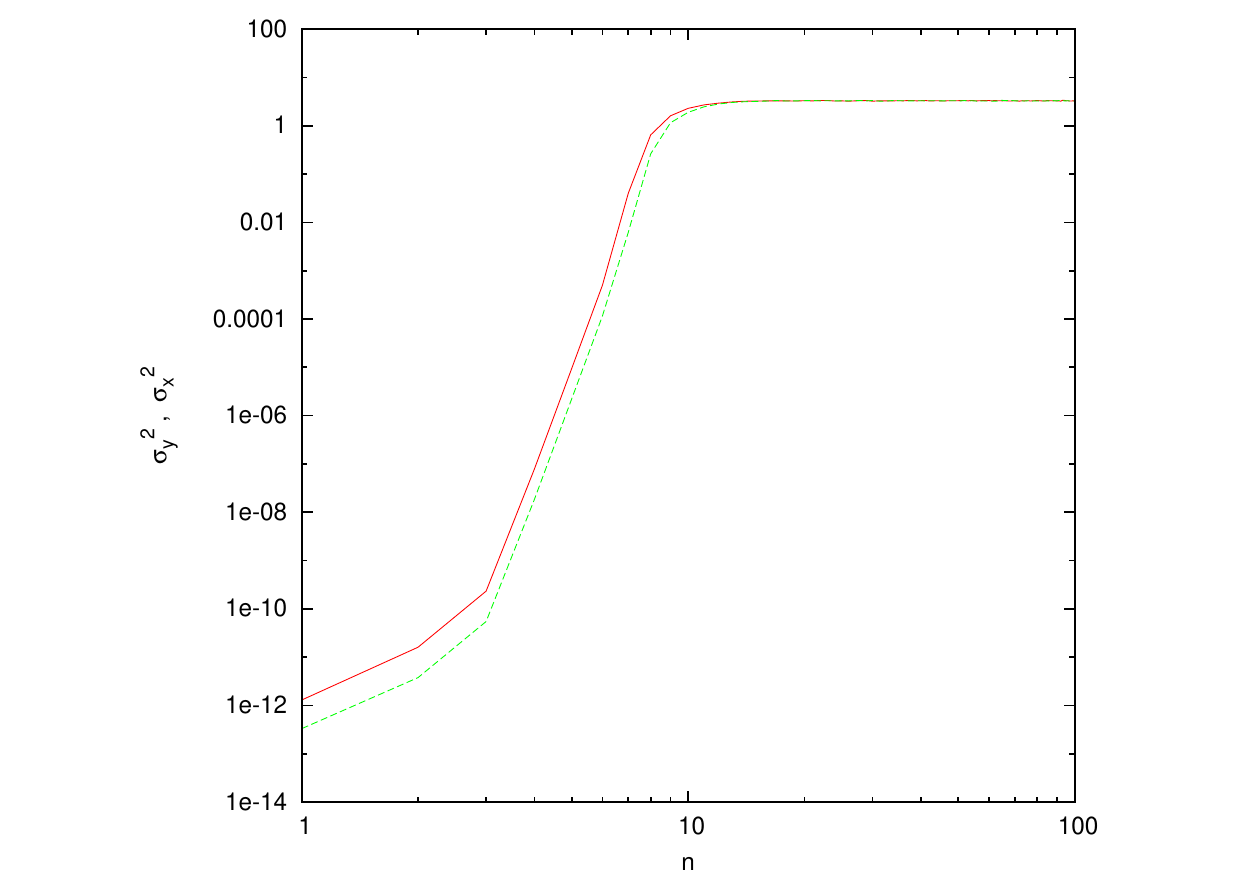} 
  \hspace{1cm}
  \includegraphics[width=0.45\textwidth]{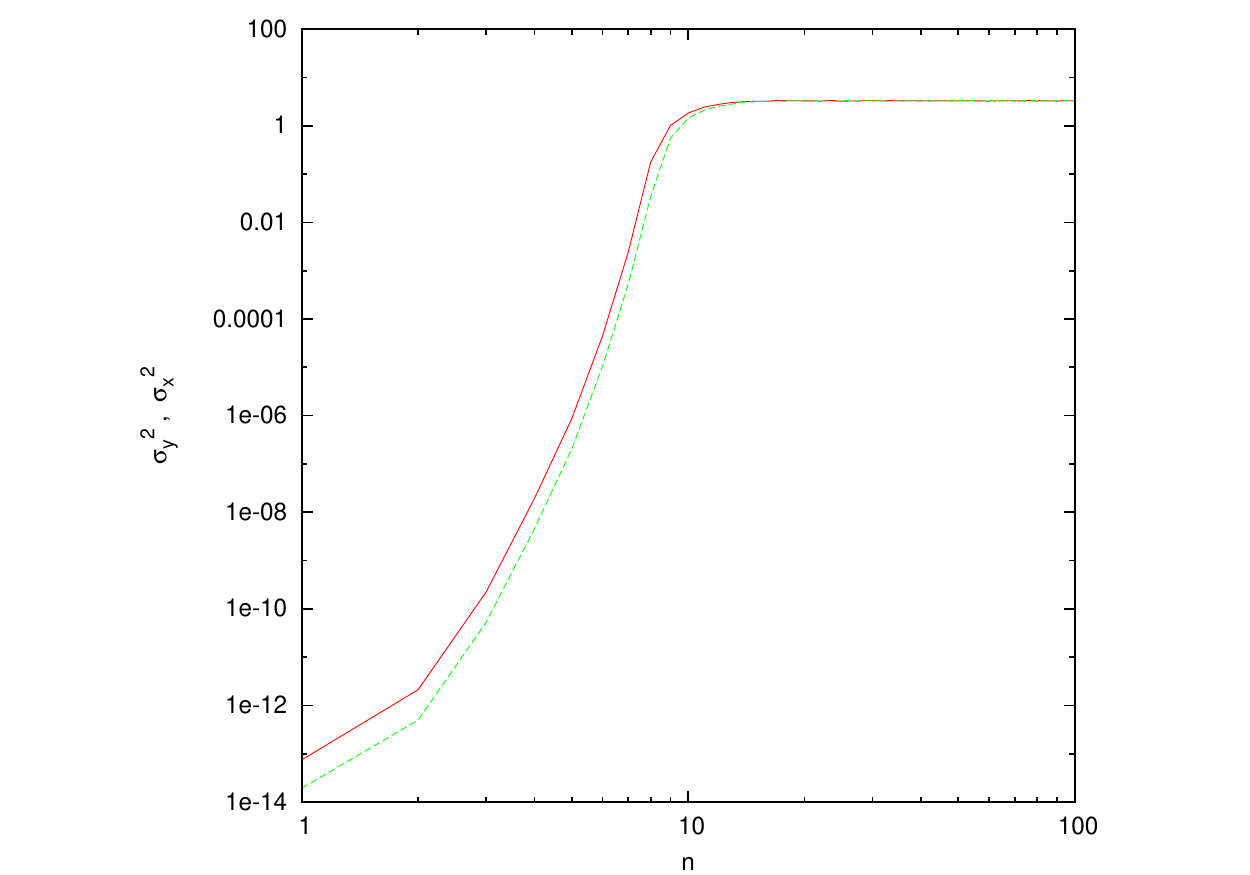}
  \caption{Evolution of the variance of a distribution of reversibility errors in action (red) and angle (green) 
    variables for the standard map with $\lambda=10$. There is evident equivalence between  round off error 
    (left) and stochastic uncorrelated perturbations of amplitude $10^{-7}$ (right).}
  \label{fig: rever_chaotic}
\end{figure}

In the fully chaotic regime for the standard map  $\lambda \gg1 $ the  random perturbation 
produces very similar results to the round off error as shown in Figure~\ref{fig: rever_chaotic}
for $\lambda=10$ since in presence of a chaotic dynamics the round off error correlations are lost.



The results for the standard map in the small $\lambda$ regime can be compared to the variances for the skew map:
\begin{equation} 
  \label{skew_map}
  y_{n+1}=y_n \mod 1; \qquad   
  x_{n+1}= x_n+y_{n+1} \mod 1. 
\end{equation}
to which the standard map in equation~(\ref{standard_map}) reduces for
$\lambda=0$.  If the skew map is randomly perturbed:
\begin{equation} 
  \label{skew_map_random}
  y_{n+1}=y_n+ \epsilon \xi_n \mod 1; \qquad   
  x_{n+1}= x_n+y_{n+1}+\epsilon \chi_n \mod 1,
\end{equation}
where $\xi_n$ and $\chi_n$ are random uncorrelated variables the growth of the variances $\sigma_y^2$ and    $\sigma_x^2$
follows a linear and cubic law respectively \citep{turchetti2006hamiltonian}. The random perturbation of a standard map with a small value of $\lambda$ (shown in Figure~\ref{fig: rever_variances}-right) shows exactly the same growth.

In the case of the round off error the behavior of the skew map with respect to the standard map with a small value of $\lambda$ 
is  quite different. Indeed, the round off error affects the  skew map just as the translation 
on the 1-D Torus: it was observed that the global error grows linearly and the variance saturates at a very small 
value with respect to the size of the torus (see Figure 7 in \citet{turchetti2010asymptotic}).
For the standard map, the coupling between action and angle, even for very small $\lambda$, causes a growth of the variance of the fluctuations
due to round off as shown in Figure ~\ref{fig: rever_variances}-left. As a consequence the effect of round off in a  very weakly perturbed standard map is similar to a random perturbation.

\section{A 4D Map}

\begin{figure}[ht!]
  \centering  
  \includegraphics[width=0.45\textwidth]{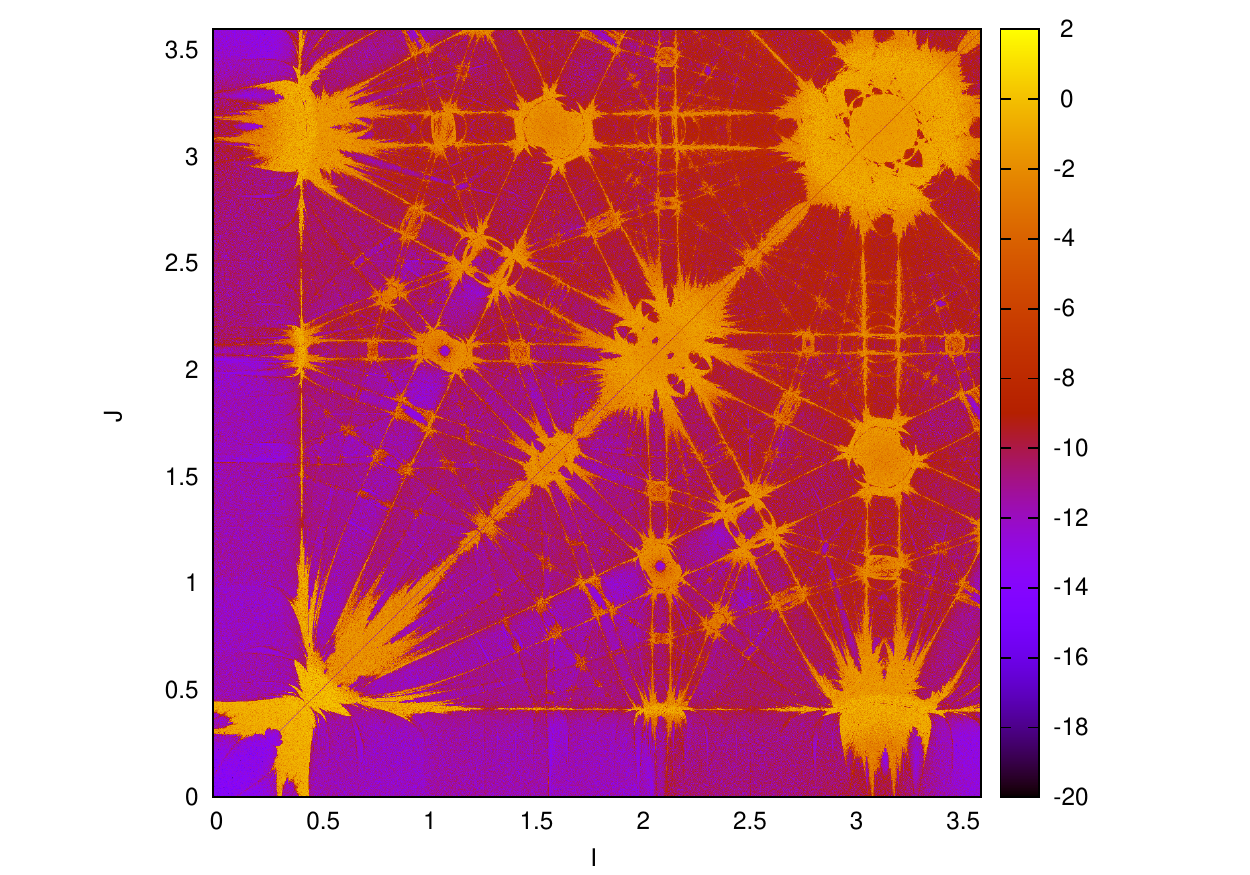}
  \hspace{1cm}
  \includegraphics[width=0.45\textwidth]{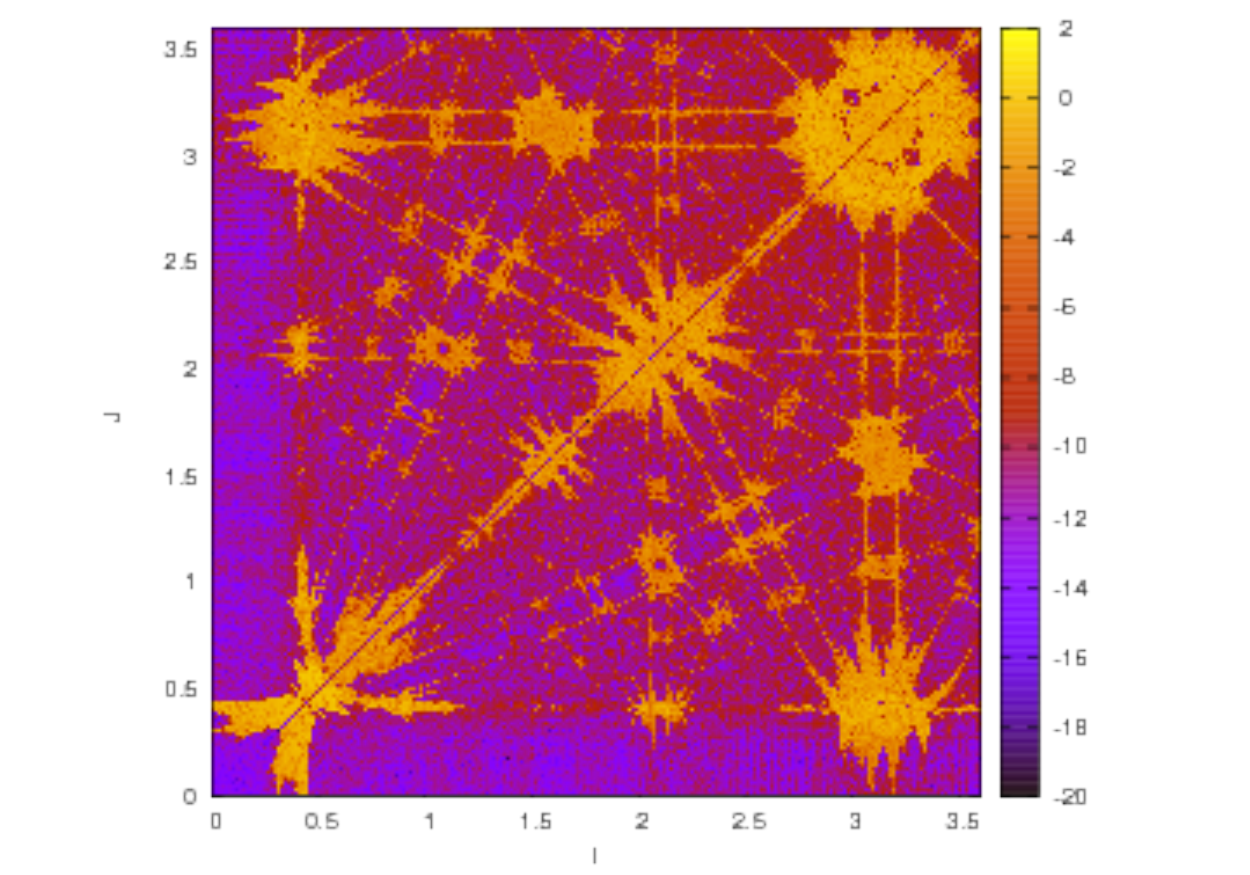}
  \caption{$\ln(R_{10^3})$ at left and $\ln(\Delta_{10^3})$ at right for map (\ref{Froeschle_map}) using $c=2$ and $\mu=0.6$. }    
  \label{fig: froes_map}
\end{figure}

\begin{figure}[ht!]
  \centering  
  \includegraphics[width=0.45\textwidth]{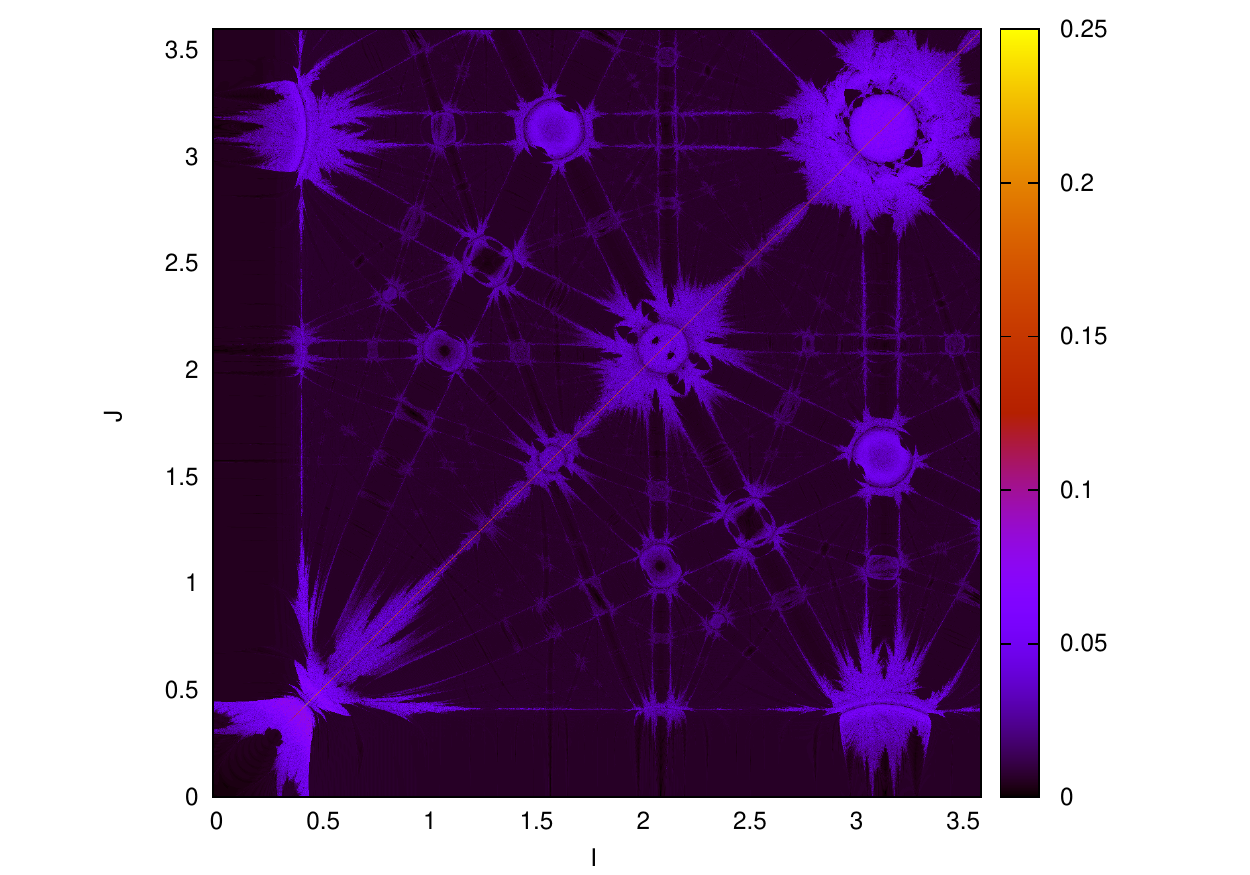}
  \hspace{1cm}
  \includegraphics[width=0.45\textwidth]{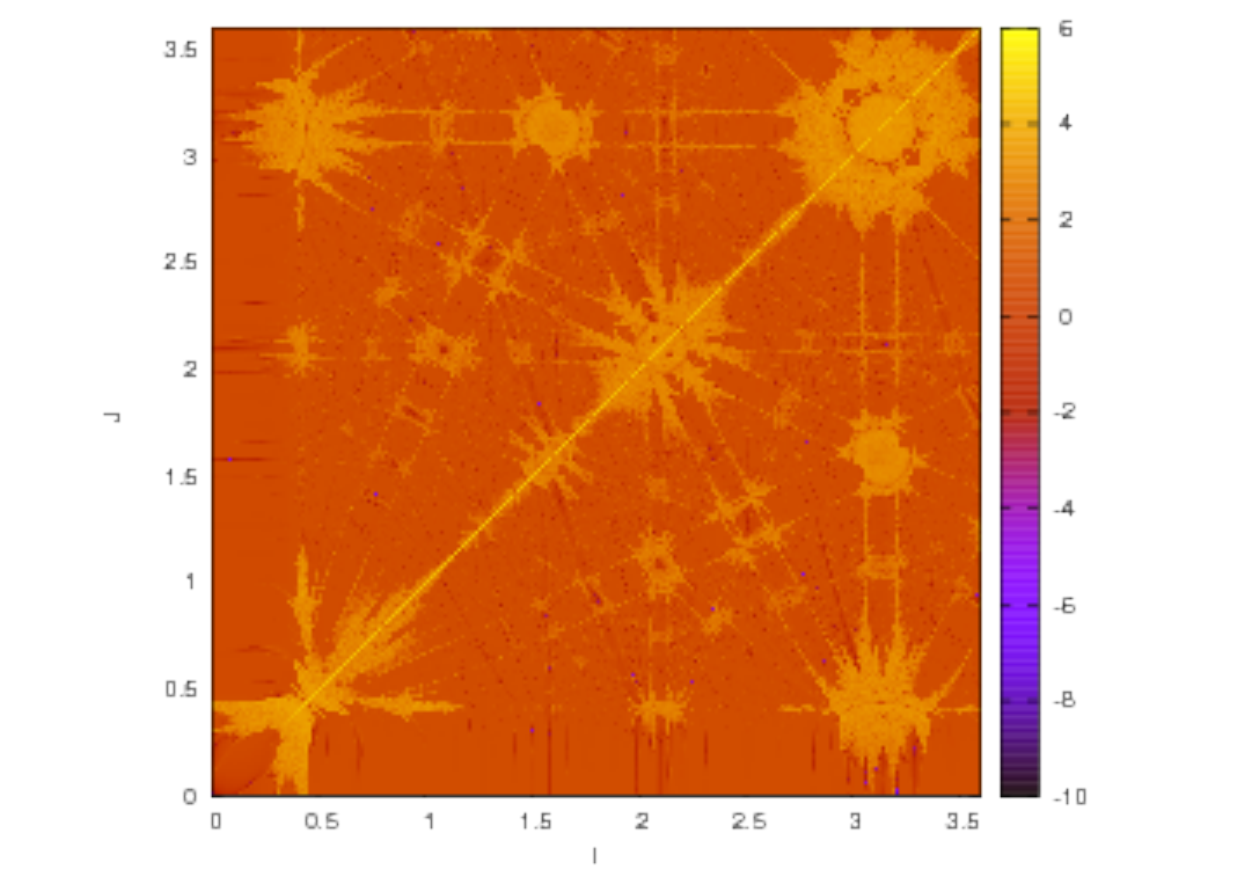}
  \caption{$\mbox{mLCE}(10^3)$ at left and $\ln(\bar{Y}(10^3))$ at right for map (\ref{Froeschle_map}) using $c=2$ and $\mu=0.6$. }    
  \label{fig: froes2_map}
\end{figure}

\begin{figure}[ht!]
  \centering  
  \includegraphics[width=0.45\textwidth]{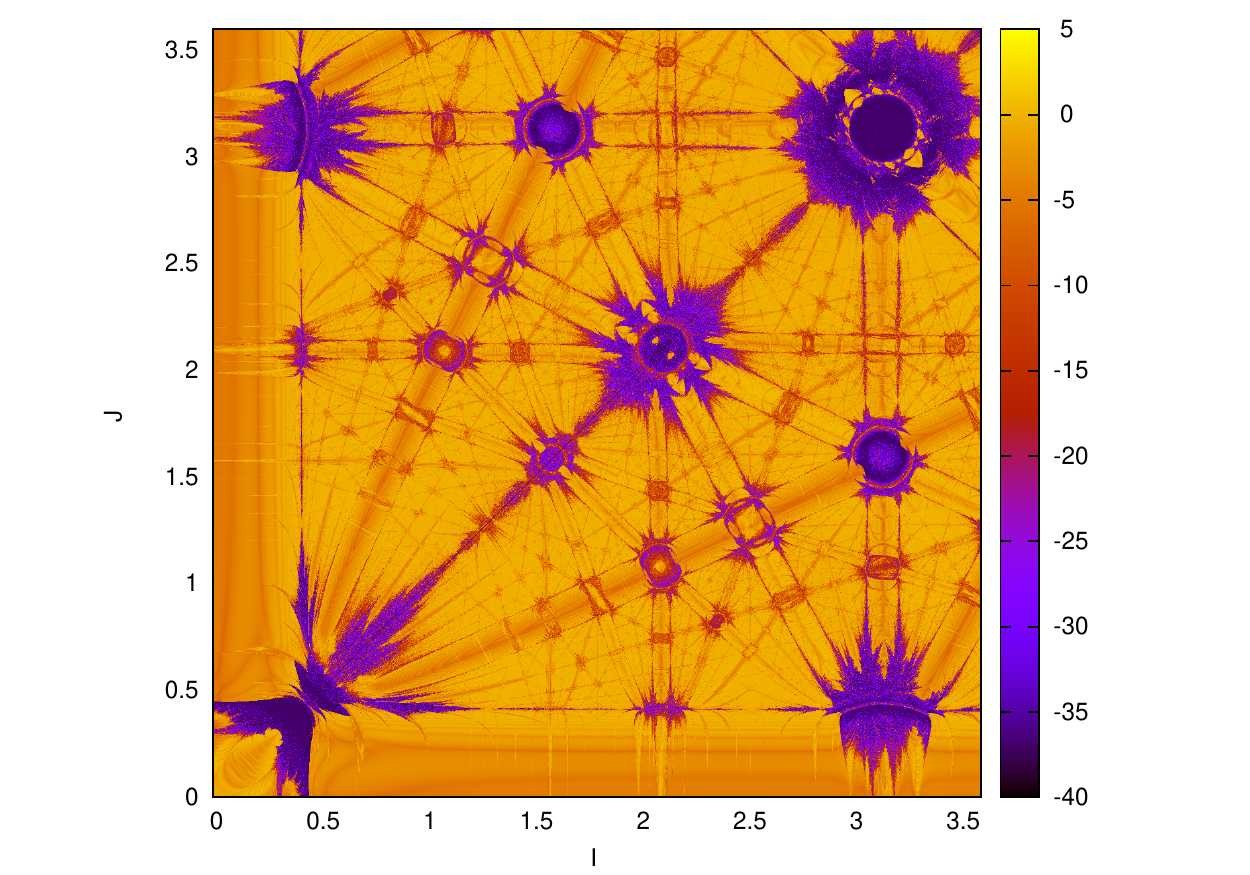}
  \caption{$\mbox{SALI}(10^3)$ for map (\ref{Froeschle_map}) using $c=2$ and $\mu=0.6$. }    
  \label{fig: froes3_map}
\end{figure}

In this section we show how either the reversibility error or the divergence of orbits can be used to 
analyze the resonance structure of four dimensional non integrable maps. As an example we consider a symplectic nearly integrable
map extensively used in the literature (see \citet{guzzo2004first}). This map  is defined as:

\begin{align}
  \label{Froeschle_map}
  \theta_{n+1} &= \theta_n + I_n \nonumber\\
  \phi_{n+1}   &= \phi_n + J_n  \nonumber\\
  I_{n+1}      &= I_n - \mu \frac{\partial V(\theta _{n+1},\phi_{n+1}) }{\partial \theta_{n+1}}\nonumber\\
  J_{n+1}      &= J_n - \mu \frac{\partial V(\theta _{n+1}, \phi_{n+1})}{\partial \phi_{n+1}};
\end{align}

where $V\equiv 1/( \cos (\theta)+\cos (\phi)+2+c )$, with $c>0$ and $\mu$ is the perturbation parameter.
We have taken a grid of $1146 \times 1146$ initial conditions  with actions $(I,J)\in[0,~3.6]\times[0,~3.6]$ 
and a fixed pair of angles, namely $(\theta,\phi)=(0.5,~0.5)$ and computed $R_n$, $\Delta_n$, mLCE(n), $\bar{Y}(n)$ and 
$\mbox{SALI}(n)$ for $n=10^3$. 
Associating these values to each initial condition while using a chromatic scale we have performed 
figures ~\ref{fig: froes_map}, ~\ref{fig: froes2_map} and  ~\ref{fig: froes3_map}, where the parameter values
$c=2$ and $\mu=0.6$ have been used. 
As we did in the previous section, in all except for the mLCE we have used the natural logarithm of the 
absolute value of the concomitant indicator and again we have used a cut off value at SALI $= 10^{-16}$.
Figure ~\ref{fig: froes_map} was done taking into account the Euclidean error in only the action plane.
As it happened for the standard map both reversibility error and orbit divergence present the same order 
of magnitude. The resonance web appears in every figure and its structure is the same. From a qualitative
viewpoint no substantial differences are found and one might conclude that the dynamical information extracted
from the reversibility error is the same as for the other four dynamical indicators we have considered. In the case of mLCE, the highest values are located on the
diagonal $I = J$ even if they are not well visible in the plot. This is anyway coherent with the fact that on
the diagonal we find the highest instability.\\
The frequency map analysis was not presented for comparison because it is computationally heavier. The error reaches its highest values in small neighbourhoods of single resonance lines, because it detects the perturbed separatrices of the
pendulum models that one could associate with the resonances, and in relatively large neighbourhoods of the intersections of resonances, where the dynamics are widely non integrable, as shown by the computation
of the interpolating Hamiltonian. 
These results are certainly not exhaustive but show that the behavior of the reversibility error is strictly related 
to the divergence of orbits and consequently it is very weak in the integrable regions where it does not have a diffusive
character as for a random error.


\section{Conclusions}

We have examined the orbit divergence and the reversibility error in 
order to determine the effect of round off error for invertible maps. 
The knowledge of the exact map is not required to compute the reversibility error and the results 
obtained are about the same with respect to the case in which we compute the orbit divergence.

By choosing an ensemble of initial data we have examined the statistics of the fluctuations due 
to round off with respect to the time average of the  error. There are different behaviours according
to the degree of chaoticity that characterises the ensemble.
For chaotic orbits the variances have an exponential growth similar to the one observed for random 
perturbations and correspondingly the decay of fidelity is super-exponential in both cases. 
For regular or quasi-regular orbits differences are observed between the effects generated by 
round off errors and random perturbations.
For a quasi-integrable map the random perturbations produce a growth of variances linear and cubic  
for actions and angles, respectively. On the other hand, the round off errors produce an initial transient 
that possibly lasts the time that the round off errors of the ensemble orbits need to decorrelate. 
After this initial transient, the behaviour of variances affected by round off and noise are similar to each other.
Finally, for a regular map such as the translation of the 1D torus, the round off variance, after a linear growth,
 saturates before the distribution of errors fills all the torus, unlike in the randomly perturbed case 
where the growth is always linear until saturation of the full torus.
The fidelity has a power law decay for round off whereas it decays exponentially for random perturbations.
This is a clear signature of the correlation between the errors due to round off.

To conclude, the reversibility error provides basically the same information as the divergence of 
orbits and it is easily accessible from a computational view point. 
This is due to the facts that it does not require the solution of the variational equations and that 
both the forward and backward orbits can  be computed using single precision.
When the initial point is varied and the number of iterations is kept fixed the reversibility error 
due to round off provides an insight on the dynamical structure of the map.
For the standard map the reversibility error provides a picture of the dynamical behavior of 
the map where not only large scale features but also small scale details can be detected. 
In the case of a 4D symplectic map the reversibility error in action space provides a similar 
picture where the resonance web and the nearby regions of weakly chaotic motions can be easily highlighted.

Even if no really new information with respect to the standard indicators is provided by the reversibility error,  we point out that this type of analysis
takes into account not only the dynamics of the map but also the unavoidable effect of finite accuracy due to round off. Moreover, by increasing
the accuracy of numerical computations, the time on which the reversibility test is computed can be increased and finer details on the phase 
space structure can be detected.

\section{Aknowledgements}

We are gratefull to the anonymous reviewers for their important comments.
The stay of MM at the Physics Department of the University of Bologna was fully 
supported by a grant from the Erasmus Mundus External Cooperation Window Lot 16 Programme, EADIC, 
financed by the European Commission.  Besides, MM acknowledges financial support of a grant 
from the Consejo Nacional de Investigaciones Cient\'ificas y T\'ecnicas de la Rep\'ublica Argentina 
(CONICET). Also, he is grateful to A. Bazzani for his helpfull discussions. DF  acknowledges the financial support of the EU FP7-ERC project NAMASTE Thermodynamics of the Climate System.

\bibliographystyle{ws-ijbc}  
\bibliography{salonicco}

\end{document}